\documentclass[11pt]{article}
\usepackage{amsmath,amssymb,mathtools}
\newtheorem{theorem}{Theorem}
\newtheorem{proposition}[theorem]{Proposition}
\newtheorem{corollary}[theorem]{Corollary}
\newtheorem{lemma}[theorem]{Lemma}
\newtheorem{example}[theorem]{Example}
\newtheorem{remark}[theorem]{Remark}
\newtheorem{definition}[theorem]{Definition}
\newtheorem{algorithm}[theorem]{Algorithm}
\usepackage{xcolor}

\newcommand{\ba}{\begin{array}}
\newcommand{\ea}{\end{array}}

\newcommand{\dm}[1]{ {\displaystyle{#1} } }

\def\R{{\mathbb R}}
\def\C{{\mathbb C}}

\newcommand {\proof} {\par{\it Proof}. \ignorespaces}
\newcommand {\eproof}
      {\space
        {\ \vbox{\hrule\hbox{\vrule height1.3ex\hskip0.8ex\vrule}\hrule}}
        \par}

\def\lam{\lambda}
\def\sig{\sigma}
\def\Lam{\Lambda}

\def\diag{\mathrm{diag}}

\def\rank{\mathrm{rank}}

\def\rar{\rightarrow}
\def \rev{\mathrm{rev}}
\def \cA{\mathcal{A}}

\newcommand {\mat}  [1] {\left[\begin{array}{#1}}
\newcommand {\rix}      {\end{array}\right]}

\title{Structure-preserving deflation of critical eigenvalues in quadratic eigenvalue problems associated with damped mass-spring systems
}

\author{Rafikul Alam \thanks{Department of Mathematics, IIT Guwahati, Guwahati-781039, India. ({\tt rafik@iitg.ac.in, rafikul68@gmail.com}) orcid: 0000-0002-1982-1955} \and
Volker Mehrmann \thanks{Institut f\"{u}r Mathematik, TU Berlin, Str. d. 17. Juni 136, D-10623 Berlin, Germany. ({\tt mehrmann@math.tu-berlin.de}) orcid: 0000-0001-5051-2870  }
\and Ninoslav Truhar \thanks{Josip Juraj Strossmayer University of Osijek, School of Applied Mathematics and Informatics, Trg Ljudevita Gaja 6, HR-31000 Osijek, Croatia,
  ({\tt ntruhar@mathos.hr}) orcid:0000-0002-0848-8309}
}

\begin{document}
\maketitle

\begin{abstract} For a quadratic matrix polynomial
associated with a damped mass-spring system there are three types of critical eigenvalues, the eigenvalues $\infty$ and $0$ and the eigenvalues on the imaginary axis. All these are on the boundary of the set of (robustly) stable eigenvalues. For numerical methods, but also for (robust) stability analysis, it is desirable to deflate such eigenvalues by projecting the matrix polynomial 
to a lower dimensional subspace before computing the other eigenvalues and eigenvectors.  We describe  structure-preserving deflation strategies that deflate these eigenvalues
via a trimmed structure-preserving linearization.
We employ these results for the special case of hyperbolic problems.
We also analyze the effect of a (possibly low rank) parametric damping matrix  on purely imaginary eigenvalues.
\end{abstract}



\renewcommand{\thefootnote}{\arabic{footnote}}
\setcounter{footnote}{0}

\section{Introduction} \label{Sec:Introduction}
We consider the quadratic eigenvalue problem
\begin{equation} \label{qep}
P(\lam)v = 0,
\end{equation}
for regular matrix polynomials  $P(\lam) := \lam^2 M+ \lam D + K$, with Hermitian coefficients $M=M^*,\ D=D^*,\ K=K^* \in \mathbb C^{n\times n}$, eigenvalues $\lambda$ and right eigenvectors $v\in  \mathbb C^{n}$. (Here $M^*$ denotes the conjugate transpose of the matrix $M$.)
Such quadratic eigenvalue problems are well studied \cite{GohLR82,Lan02,Prz68}
and arise in a multitude of applications; see~\cite{BetHMST13,MehV05,TisM01}.
Our main motivation to revisit this problem is the class of damped mass-spring systems, see \cite{Ves11}, where  $M$ is a positive definite or positive semidefinite  mass matrix, $D$ is a semidefinite damping matrix and  $K$ is  a positive (semi-)definite stiffness matrix.
Most of our  results hold regardless of whether $D$ is positive or negative semidefinite. In this paper we focus on the case that $D$ is positive semidefinite.

Homogeneous damped mass-spring systems have the structure
\begin{equation}\label{sos}
    M \ddot x + D \dot x + K x =0,
\end{equation}
and the ansatz $x=e^{t\lambda} v$ leads to the quadratic eigenvalue problem \eqref{qep}.
Numerical methods for the solution of the eigenvalue problem are widely available;
see~\cite{GueT17,MehV05}. They are mainly based on the reformulation (\emph{linearization}) of \eqref{qep} as a linear eigenvalue problem.
The construction of such linearizations has been a very important research topic (see, e.g.,
\cite{DopLPV18,MacMMM06a}) in particular when it is essential to preserve the  structure, see~\cite{HigMMT07,MacMMM06b}.

One of the difficulties with such structured linearizations arises when the leading coefficient $M$ and/or the trailing coefficient $K$ is singular, which may mean that the classical structured linearizations do not exist \cite{MacMMM06b}, and the occurring eigenvalues at $\infty$ or $0$ may have Jordan blocks of size $2$, or the problem is non-regular, see \cite{MehMW18,MehMW22a}.

In such a situation it is advisable to first deflate the part of the matrix polynomial associated with these eigenvalues. This can, for example,  be done via \emph{trimmed} linearizations  \cite{ByeMX08}. However, it is common practice in industrial applications to introduce small perturbations  to move the eigenvalues away from $0$, $\infty$ but this may lead to drastically wrong results, see   the analysis in the context of brake squeal \cite{GraMQSW16}. It has been shown there, that from a numerical point of view it is better to consider the nearby problem with exact eigenvalues $0$ or $\infty$ and to apply numerical methods to a problem where these eigenvalues have been deflated.
Let us illustrate the situation of eigenvalues near $\infty$ or $0$ with a very simple example, see \cite{MehS23}.

\begin{example} \label{exmsd}{\rm
  Consider the standard model of a damped mass-spring system
$ m\ddot q + d \dot q + kq=f$,
with position mass $m$, stiffness $k$ and damping $d$, forcing function $f$ and kinetic plus potential energy $\mathcal H(q,p) = \frac{1}{2} kq^2 + \frac{p^2}{2m}$. The classical first order formulation (by introducing $p=\dot q$) is given by
\[
\mat{cc}  m & 0 \\ 0 & 1\rix\mat{c} \dot{p} \\ \dot{q} \rix = \mat{cc}  -d & -k \\ 1 & 0\rix  \mat{c} p \\ q \rix + \mat{c} f \\ 0 \rix.
\]
%
%
The associated matrix pencil has the
%
eigenvalues
$ -\frac{d}{2m}\pm\sqrt{(\frac{d}{2m})^2-\frac{k}{m}}$.

If the size of the mass $m$ is negligible compared to the damping and stiffness, it is common practice to set $m=0$, so that the associated matrix pencil
%
has eigenvalues $\infty$ and $-k/d$ or even a Jordan block at $\infty$ if $d=0$.

Another common technique is to replace stiff springs (with very large $k$) by rigid connections. In our example this corresponds to the limit  $k \to \infty$. To be able to take this limit we can rewrite the system in new coordinates and obtain an equivalent matrix pencil with
\[
P(\lam) = \lam \mat{cc} m & 0 \\ 0& \frac{1}{k}\rix +\mat{cc} d & 1 \\ -1
& 0 \rix.
\]
For $k \to \infty$ this pencil has a Jordan block of size $2$ at the eigenvalue $\infty$.
One may also take the limits $m \to 0$ and $k \to \infty$ simultaneously. This gives a double semisimple eigenvalue at $\infty$.
}
\end{example}

Consider the linearization of \eqref{qep} as
\begin{equation}\label{ph} Q(\lam) := \lam \left[\begin{matrix}
M &  \\ &  I_n
\end{matrix} \right] + \left[ \begin{matrix}
 D & I_n \\ -I_n & 0
\end{matrix}\right]\left[\begin{matrix}
I_n &  \\ &  K
\end{matrix} \right]
,
\end{equation}
with $I_n$ denoting the identity matrix of size $n$. This has the structure of a \emph{dissipative Hamiltonian pencil} and the spectral properties for such pencils have been
analyzed in detail in
\cite{MehMW18}.
If $M$ (or $K$) is invertible, the case that we are considering in this paper, then it has been shown there that the Jordan blocks at $0$ (resp., $\infty$) are of size at most one. Further properties of the dynamical system are studied in \cite{AchAM21,MehMW22a,MehU23}.

In view of the possible Jordan structures that may arise when the matrix polynomial has eigenvalues at $0$ and/or $\infty$, one of the goals of this paper is to derive { structure preserving} deflation methods for these eigenvalues. This is done for general damped mass-spring systems in Section~\ref{sec:deflation} and for the specific case of \emph{hyperbolic problems}  in Section~\ref{sec:hyperbolic}. { The main reason for preserving the structure is that it leads to well-known spectral properties, see e.g. \cite{MehMW18,MehMW22a},  that would get destroyed if unstructured deflation procedures are used, see e.g. \cite{MehX14a} for an extensive discussion and  perturbation analysis.}

Besides the eigenvalues $0$ and $\infty$, we also study  the  nonzero purely imaginary  eigenvalues which correspond to undamped oscillatory solutions.  We discuss deflation procedures for these  eigenvalues in Section~\ref{Sec:Evolutionofpurimageig} and analyze the effect of (low rank) damping on these eigenvalues. We also develop a constructive method for low rank perturbations of the damping matrix  that move all or a selected set of purely imaginary eigenvalues of $P(\lam)$ from the imaginary axis.

We denote the set of all $m\times n$ matrices with entries in $\mathbb{C}$ by $\mathbb{C}^{m\times n}.$
The nullspace of a matrix $X$ is denoted by $\mathcal N(X)$. The dimension of the nullspace is called the \emph{nullity} of $X$.
If $ X$ and $Y$ are $n\times n$ Hermitian matrices, then we write $ X> Y$ (resp., $X\geq Y$) if $ X-Y$ is positive definite (resp., positive semidefinite).
The \emph{inertia of a Hermitian matrix} $ H$ is  denoted by $ \imath(H) := ( \imath_+(H), \, \imath_{-}(H), \,\imath_0(H))$, where $\imath_+(H), \, \imath_{-}(H)$ and $ \imath_0(H)$ denote the number of positive, negative and zero eigenvalues of $H$, respectively.

We assume that the matrix polynomial $P(\lam)$ in (\ref{qep}) is \emph{regular}, i.e., there exists $\lambda_0  \in \mathbb C$ such that $\det  P(\lam_0)\neq 0$. The \emph{spectrum} of $P(\lam)$, denoted by $\Lam(P),$ is given by
$ \Lam(P) :=\{ \lam \in \C : \det(P(\lam)) =0\}
$.
Let
\[
\mathbf{S}_P(\lam) = \diag( \phi_1(\lam), \ldots, \phi_n(\lam))
\]
be the \emph{Smith form}, see, e.g.,~\cite{Ros70}, of $P(\lam)$ which is obtained under unimodular equivalence transformations, where $ \phi_1, \ldots, \phi_n$ are unique monic polynomials and $\phi_i $ divides $\phi_{i+1}$ for $ i =1,\ldots, n-1$.
Then $\phi_{i}(\lam) = (\lam - \mu)^{m_i}\rho_{i}(\lam)$ with $\rho_{i}(\mu) \neq 0$ and $m_i \geq 0$ for $i = 1,\ldots,n$.  The tuple
$\textrm{Ind}(\mu, P) := (m_1, \ldots, m_n)$ is called the {\em multiplicity index } of $\mu$ and satisfies the condition $0\leq m_1 \leq m_2 \leq \cdots \leq m_n. $ The nonzero components in $\textrm{Ind}(\mu, P)$ are called the {\em partial multiplicities} of $\mu$ as an eigenvalue of $P(\lam)$. The factors $(\lam - \mu)^{m_i}$ with $m_i \neq 0$ are called the {\em elementary divisors} of $P(\lam)$ at $\mu$. The {\em algebraic multiplicity} of $\mu$ is then given by $m_1+\cdots+m_n$.

\begin{definition}\label{trim} Let $ P(\lam) := \lam^2 M + \lambda D+ K$ be an $n\times n$ matrix polynomial, where $M$ is nonsingular and $\rank(K) < n.$ An $m\times m$ matrix pencil $L(\lam) :=  {L_0} + \lambda {L_1}$, where $L_1$ is nonsingular, is said to be a trimmed linearization of $ P(\lambda)$ if    $m < 2n, \Lambda(P) \setminus \{0\} = \Lambda(L)\setminus\{0\}$ and $\textrm{Ind}(\mu, P) = \textrm{Ind}(\mu, L)$ for all $ \mu \in \Lambda(P) \setminus\{0\}.$
\end{definition}

For $ P(\lam) := \lam^2 M + \lambda D+ K$ with $\rank(M) < n$ and $K$ nonsingular,  a trimmed linearization of $P(\lam)$ is defined similarly. 

\begin{remark}\label{rem:trimmed}{\rm   This definition of a trimmed linearization differs from the original definition in \cite{ByeMX08}, where all the chains of length larger than $1$ associated with the eigenvalue $0$ are deflated. Here we allow any deflation of eigenvalues that lead to a size reduction in the linearization. However, the trimmed linearizations in [6] are a special case of the trimmed linearizations that we consider here.}
\end{remark}

Let $G(\lam)$ be an $n\times n$ regular rational matrix. Then the Smith-McMillan form~\cite{Ros70} of $G(\lam)$, which is obtained under unimodular equivalence transformations, is given by
\[
\mathbf{SM}(\lam) := \diag( \phi_1(\lam)/{\psi_1(\lam)}, \ldots, \phi_n(\lam)/{\psi_n(\lam)}),
\]
where $ \phi_1, \ldots, \phi_n$ and $\psi_1, \ldots, \psi_n$ are unique monic polynomials such that $\phi_i$ and $\psi_i$ are pairwise coprime, $\phi_i$ divides $\phi_{i+1}$ and $ \psi_{i+1}$ divides $\psi_i.$

A matrix polynomial  $P(\lam)$ as in \eqref{qep} is said to be \emph{hyperbolic} if  $ M>0$ and
\begin{equation}\label{defhyp}
(x^*Dx)^2 -  4 x^*Mx \, x^* Kx > 0
\end{equation}
for all nonzero $x$, or equivalently, if $ - P(\mu) > 0$ for some  $ \mu \in \mathbb{R}$,  see~\cite{HigMT09}. If $P(\lam)$ is hyperbolic with $D >0$ and $ K\geq 0$ then $P(\lam)$ is called \emph{overdamped.}

An $n \times n$ Hermitian matrix pencil $A+\lam E$ is said to be a \emph{definite pencil} if $ A+\mu E >0$ for some $ \mu \in \R\cup\{\infty\}$.
The \emph{reverse polynomial} to  $P(\lam)$ in \eqref{qep} is $\rev P(\lam):=\lam^2 K+\lam D+M$.
For $M=M^*>0$, a matrix $X \in \C^{n\times n}$ is called \emph{$M$-unitary} if $ X^*MX = I_n$.

\section{Deflation of zero and infinite eigenvalues}\label{sec:deflation}

In this section we design a method for deflating zero (infinite) eigenvalues of $P(\lam)$ and we assume that either $M$ or $K$ is invertible.    Our strategy is to construct \emph{trimmed} structured linearizations of $P(\lam)$ in which the zero or infinite eigenvalues of $P(\lam)$ have been deflated. We treat only the quadratic matrix polynomial $P$ but the proposed method can be extended to general even order Hermitian matrix polynomials.

Since $\infty$ is an eigenvalue of $P(\lam)$ if and only if $0$ is an eigenvalue of $\rev P(\lam):=\lam^2 K+\lam D+M$,  without loss of generality we describe the procedure for deflating the zero eigenvalues, i.e., we discuss the case that $ M >0$ and that $K \geq 0$ has rank $r < n$.
Although most of our results hold for $P(\lam)$ when $ M = M^*$ is nonsingular and $D=D^*$ is semidefinite, to not overload the presentation, we only consider  $ M >0$ and $ D\geq 0$. Furthermore,  we assume in the following that we have a full rank factorization $ K = GG^*$, i.e., $G \in \C^{n\times r}$  has full column rank.
In many applications this is directly available, \cite{Ves11}, otherwise it may be obtained, e.g.,  by a low rank Cholesky  factorization of $K$.
Then  we may  consider the  Hermitian matrix pencil
\begin{equation}\label{herm} S(\lam) := \left[ \begin{matrix}
\lam M + D & G \\ G^* & - \lam I_r
\end{matrix}\right] = \lam \left[\begin{matrix}
M &  \\ &  -I_r
\end{matrix} \right] + \left[ \begin{matrix}
 D & G \\ G^* & 0
\end{matrix}\right]
\end{equation}
of size $n+r$ where both coefficients are indefinite.

\begin{remark}\label{rem:DG}{\rm
As an alternative to \eqref{herm}, one may consider the
equivalent \emph{dissipative Hamiltonian  pencil}
\begin{equation}\label{jherm} L(\lam) := \lam \left[\begin{matrix}
M &  \\ &  I_r
\end{matrix} \right] + \left[ \begin{matrix}
 D & G \\ -G^* & 0
\end{matrix}\right].
\end{equation}
%
}
\end{remark}

Our first result characterizes the Smith form of $S(\lambda)$ in relation to the Smith-McMillan form of the rational matrix function $T(\lambda):=P(\lambda)/\lambda$.
\begin{lemma}\label{lem:lem1}  Let $ P(\lam) := \lam^2M+\lam D+ GG^*$, where $0<M=M^*,\ 0 \leq D=D^*\in \mathbb C^{n\times n}$, and $ G \in \mathbb C^{n\times r}$ is of rank $r$. If  the \emph{Smith-McMillan form} of the rational matrix  $T(\lam)$ is
$\diag( \phi_1(\lam)/{\psi_1(\lam)}, \cdots, \phi_n(\lam)/{\psi_n(\lam)})$,
%
then the \emph{Smith form} of the matrix pencil $S(\lam)$ in \eqref{herm} is given by
\begin{equation}\label{smform}
\diag(I_r,\phi_1(\lam), \ldots, \phi_n(\lam))
\end{equation}
In particular,
 $ \Lam(S)\setminus\{0\} = \Lam(P)\setminus\{0\}$  and $\mathrm{Ind}(\mu, S) = \mathrm{Ind}(\mu, P)$ for all $\mu \in \Lam(P)\setminus\{0\}$.

 For any $\lambda_0\in \mathbb C\setminus\{0\}$, the map
 \[
 \mathcal N(P(\lam_0)) \longrightarrow \mathcal N(S(\lam_0)),\ u \longmapsto \left[\begin{matrix}
\lam_0 u \\ G^*u
\end{matrix}\right],
\]
is an isomorphism.

Furthermore, if  $D$ is invertible, then $ 0 \in \Lam(S)$ if and only if $ G^*D^{-1}G$ is singular. In particular, if $D$ is definite then $ \Lam(S) = \Lam(P)\setminus\{0\}$.
\end{lemma}
\proof We have
\begin{eqnarray*} \det(S(\lam)) &=& \det( - \lam I_r) \det( \lam M +D+ G(\lam I_r)^{-1}G^*)  \\ &=&  (-1)^r \lam^r \det(P(\lam)/\lam) = (-1)^r \det(P(\lam))/{\lam^{n-r}},
\end{eqnarray*}
which shows that $\Lam(S)\setminus\{0\} = \Lam(P)\setminus\{0\}$. If $\lam_0 \in \Lam(P)\setminus \{0\}$, then it follows that the algebraic multiplicity of $\lam_0 $ as an eigenvalue of $S(\lam)$ is the same as the algebraic multiplicity of $\lam_0$ as an eigenvalue of $P(\lam)$. It is obvious that the map
\[
\mathcal N(P(\lam_0)) \longrightarrow \mathcal N(S(\lam_0)),  \quad u \longmapsto \left[\begin{matrix}
\lam_0 u \\ G^*u
\end{matrix}\right]
\]
is an isomorphism. Hence the geometric multiplicity of $\lam_0$ as an eigenvalue of $S(\lam)$ is the same as the geometric multiplicity of $\lam_0$ as an eigenvalue of $P(\lam)$.

Since
\[
\rank\left( \left[\begin{array}{cc} G^* & - \lam I_r\end{array}\right]\right) = r = \rank\left( \left[\begin{array}{c} G \\ - \lam I_r\end{array}\right]\right)
\]
for all $ \lam \in \C$,
it follows that
$
T(\lam)= \lam M +D+ G ( \lam I_r - 0)^{-1} G^*
$
is a \emph{minimal realization}, i.e., a realization of minimal degree (see, e.g., \cite{Ant05})  of the rational matrix  $T(\lam)$.
As $T(\lam)$ is the transfer function associated with the system matrix $S(\lam)$ and has Smith-McMillan form $
\diag( \phi_1(\lam)/{\psi_1(\lam)}, \ldots, \phi_n(\lam)/{\psi_n(\lam)})$, 
the Smith form of $S(\lam)$ is given by \eqref{smform}, see, \cite[Theorem~4.1, p.111]{Ros70}.
%
%
Hence it follows that $\mathrm{Ind}(\mu, S) = \mathrm{Ind}(\mu, P)$ for all $\mu \in \Lam(P)\setminus\{0\}$.

Next, we have  $\det(S(0)) = \det(D) \det( - G^*D^{-1}G) $ showing that $S(0)$ is singular if and only if $G^*D^{-1}G$ is singular.
Therefore, $0$ is an  eigenvalue of $S(\lam)$ if and only if $ G^*D^{-1}G$ is singular. In particular,  if  $D$ is definite then $ G^*D^{-1}G$ is nonsingular and the assertion follows.
\eproof

\begin{remark} \label{rem:remtrim}{\rm
Lemma~\ref{lem:lem1} shows that $S(\lam)$ in \eqref{herm} is a trimmed linearization of $P(\lam)$ in the sense of Definition~\ref{trim}. For  $D$ definite, we have $$\Lam(S) = \Lam(P)\setminus \{0\} \text{ and  }\mathrm{Ind}(\mu, P) = \mathrm{Ind}(\mu, S)$$ for all $\mu \in \Lam(P)\setminus\{0\}$. Hence, in this case, the linearization $S(\lam)$ deflates the zero eigenvalue of $P(\lambda)$ and  preserves the nonzero  eigenvalues of $P(\lam)$ including their partial multiplicities.
}\end{remark}

For the case that $D$ is only semidefinite
we need the following lemma.
\begin{lemma} \label{lem:null} Let $ \mathcal{D} :=\left[\begin{array}{cc|c} D_{11} & D_{12} & D_{13} \\ D_{21} & D_{22} & 0 \\ \hline  D_{31} & 0 &  0 \end{array}\right]$, 
where  $ D_{13},D_{31} \in \C^{r\times r}$ are nonsingular. For $ \mathcal{\widehat D} := \left[\begin{array}{cc} D_{12} & D_{13} \\ D_{22} & 0\end{array}\right]$ then
$\mathrm{nullity} (\mathcal{D}) = \mathrm{nullity}(D_{22}) = \mathrm{nullity}(\mathcal{\widehat D})$.

    Moreover, the maps
    \begin{eqnarray*}
    && \mathcal N(D_{22}) \longrightarrow \mathcal N(\widehat {\mathcal{D}}),\quad x \longmapsto \left[\begin{array}{c} x \\ - D_{13}^{-1}D_{12} x\end{array}\right], \\
    &&
    \mathcal N(\widehat {\mathcal{D}}) \longrightarrow N(\mathcal{D}),\quad  \left[\begin{array}{c} x \\ y \end{array}\right]\longmapsto \left[\begin{array}{c} 0 \\ x \\ y \end{array}\right]
    \end{eqnarray*}
    are isomorphisms. Furthermore,  $\rank(\mathcal{\widehat D}) = r + \rank(D_{22})$ and  $\rank(\mathcal{D}) = 2 r + \rank(D_{22})$.
    \end{lemma}
    \proof We have~~$\widehat{\mathcal{D}}\left[\begin{array}{cc} I_r & -D_{13}^{-1} D_{12}\\ 0 & I_{n-r}\end{array}\right] = \left[\begin{array}{cc} I_r & 0\\ 0 & D_{22}\end{array}\right]$. Hence $\rank(\widehat{\mathcal{D}}) = r + \rank(D_{22})$ and $\mathrm{nullity}(\widehat{\mathcal{D}})
    = \mathrm{nullity}(D_{22})$. It is easily seen that the map
    \[
    \mathcal N(D_{22}) \longrightarrow \mathcal N(\widehat{\mathcal{D}}),\quad  x \longmapsto \left[\begin{array}{c} x \\ - D_{13}^{-1}D_{12} x\end{array}\right]
    \]
    is well-defined and an isomorphism.

    By applying block Gaussian elimination   $\mathcal{D}$ can be  transformed to
    \[
    \left[\begin{array}{cc|c} 0 & 0 & D_{13} \\ 0 & D_{22} & 0 \\ \hline  D_{31} & 0 &  0 \end{array}\right],
    \]
    which shows that  $\rank(\mathcal{D}) = 2 r + \rank(D_{22})$ and $\mathrm{nullity} (\mathcal{D}) = \mathrm{nullity}(D_{22})$. Finally, the map from $\mathcal N(\widehat{\mathcal{D}}) $ to $\mathcal N(\mathcal{D})$ is easily seen to be well-defined and is an isomorphism.
    \eproof

    As an immediate consequence of Lemma~\ref{lem:null}, we  show  that the projection (restriction) of the damping matrix $D$ onto the nullspace of  $K := GG^*$ plays a vital role in the deflation of zero eigenvalues of $P(\lam)$.

    \begin{lemma} \label{lem:new1}  { Consider  $ P(\lam) := \lam^2M+\lam D+GG^*$ with $D=WW^*$ and the associated   pencil $S(\lam)$ in \eqref{herm},} where $0<M=M^*\in \mathbb C^{n\times n}$, $ G \in \mathbb C^{n\times r}$ is of rank $r$, and $ W \in \mathbb C^{n \times m}$ is of rank $m$. Let  $Q^*G = \left[ \begin{array}{c} R \\ 0\end{array}\right]$ be a QR factorization of $G$,  where $R \in \C^{r\times r}$ is nonsingular and
    $ Q :=\left[ \begin{array}{cc}Q_1 & Q_2\end{array}\right] $ is  unitary with $ Q_2 \in \C^{n\times (n-r)}$. Then   $ \Lam(S)  = \Lam(P)\setminus\{0\}$ if and only if $Q_2^*DQ_2$ is nonsingular.
    In particular,   $ \Lam(S)  = \Lam(P)\setminus\{0\}$ if and only if  $\rank(W^*Q_2) = n-r.$
    \end{lemma}
    \proof  We have
    \[ \diag(Q^*, I_r) S(\lam) \diag(Q, I_r)  = \lam \left[\begin{array}{cc|c} M_{11} & M_{12} & 0 \\ M_{21} & M_{22} & 0 \\ \hline  0 & 0 &  -I_r \end{array}\right] +  \left[\begin{array}{cc|c} D_{11} & D_{12} & R \\ D_{21} & D_{22} & 0 \\ \hline  R^* & 0 &  0 \end{array}\right],
    \]
    where $D_{22} = Q_2^*DQ_2$. By Lemma~\ref{lem:null}, $ 0 \in \Lam(S)$ if and only if $ D_{22}$ is singular. Hence $ \Lam(S)  = \Lam(P)\setminus\{0\}$ if and only if $D_{22}$ is nonsingular. Finally, note  that  $ Q_2^*WW^*Q_2$ is nonsingular if and only if $\rank( W^*Q_2) = n-r$.
    \eproof

    Lemma~\ref{lem:new1} shows that if $D$ is positive semidefinite,   to deflate the zero eigenvalues of $P(\lam)$,   the rank of $D$ must be greater than or equal to the nullity of $K=GG^*$. But this may not be sufficient as the following example shows.
    \begin{example}\label{ex_lem2.3}{\rm
    Consider the quadratic polynomial
    \[
    P(\lambda)=\lambda^2 M+\lambda D+K,
    \]
    with
    \[
    M=I_2,\qquad
    K=GG^*=\begin{bmatrix}1&0\\0&0\end{bmatrix},\qquad
    G=\begin{bmatrix}1\\0\end{bmatrix},
    \]
    Since $G$ is already upper triangular, we have $Q=I_2$ and
    \(
    Q_2=\begin{bmatrix}0\\1\end{bmatrix}.
    \)
    For \(
    D=\begin{bmatrix}0&0\\0&1\end{bmatrix}
    \), then
    \(
    Q_2^*DQ_2=[1],
    \)
    which is nonsingular. Moreover,
    \[
    \Lambda(P)=\{0,-1,\pm i\},
    \]
    while the associated pencil
    \[
    S(\lambda)=
    \begin{bmatrix}
    \lambda M+D & G\\
    G^* & -\lambda
    \end{bmatrix}
    \]
    has
    \[
    \Lambda(S)=\{-1,\pm i\}
    =\Lambda(P)\setminus\{0\},
    \]
    exactly as in Lemma~2.3.
    On the other hand, if \(
    D=\begin{bmatrix}1&0\\0&0\end{bmatrix}
    \), then \(
    Q_2^*DQ_2=[0], \)
    is singular and $0$ remains an eigenvalue of $S(\lambda)$.
    }
    \end{example}
     The next lemma characterizes the deflation of zero eigenvalues of $P(\lam)$ under the condition that $D$ is semidefinite  and $\rank\left(\left[ \begin{array}{cc} D & G\end{array}\right]\right) = n$.

    \begin{lemma} \label{lem:new2} { Consider  $ P(\lam) := \lam^2M+\lam D+GG^*$ and  the associated   pencil $S(\lam)$ in \eqref{herm},} where $0<M=M^*\in \mathbb C^{n\times n}$,  $0\leq D=D^*\in \mathbb C^{n\times n}$, $ G \in \mathbb C^{{n\times r}}$ is of rank $r$.     If $\rank\left(\left[ \begin{array}{cc} D & G\end{array}\right]\right) = n$ then $ 0 \notin \Lam(S)$ and we have $\Lam(S) = \Lam(P)\setminus\{0\}.$
    \end{lemma}

    \proof If $ \rank\left(\left[ \begin{array}{cc} D & G\end{array}\right]\right)= n$ then $ \mathcal N(D)\cap \mathcal N(G^*) = \{0\}$.  We now show that $\mathcal{D} := \left[ \begin{matrix}
     D & G \\ G^* & 0
    \end{matrix}\right]$ is nonsingular.  Consider the QR factorization $Q^*G =  \left[ \begin{array}{c} R \\ 0\end{array}\right]$, where $R \in \C^{r\times r}$ is nonsingular and $ Q := [Q_1, Q_2]$ is unitary with $Q_2 \in \C^{n\times (n-r)}$. Then $G^*Q_2 =0$ and the columns of $Q_2$ span $\mathcal N(G^*)$. Furthermore, we have
    \[
    \diag(Q^*, I_r) \mathcal{D} \diag(Q, I_r) = \left[\begin{array}{cc|c} D_{11} & D_{12} & R \\ D_{12}^* & D_{22} & 0 \\ \hline  R^* & 0 &  0 \end{array}\right],
    \]
    where $D_{22} := Q_2^*DQ_2$.  By Lemma~\ref{lem:null}, $\mathcal{D}$ is singular if and only if $D_{22}$ is singular.

    Now suppose that $D_{22} x = 0$ for some $x\neq 0$. Then  $ (Q_2x)^*DQ_2x =0$ and since $D$ is semidefinite,  we have $DQ_2x = 0$ which implies that $Q_2x \in \mathcal N(D)$. Since $Q_2 x \in  \mathcal N(G^*)$, we have that $ Q_2x \in \mathcal N(D)\cap \mathcal N(G^*) = \{0\}$ and thus $Q_2x = 0$ and hence $x= 0$, which is a contradiction. Therefore, $D_{22}$ is nonsingular and $ 0 \notin \Lam(S)$.
    \eproof

If $D$ is semidefinite and $\rank\left(\left[ \begin{array}{cc} D & G\end{array}\right]\right) = n-k$, then  $S(\lam)$ still has at least $k$ zero eigenvalues and extra work has to be performed for the deflation of the zero eigenvalues. To see this,  we combine the previous lemmas to construct a trimmed linearization of $P(\lam)$ that deflates all the zero eigenvalues of $P(\lam)$ when $ M>0$ and $D$ is semidefinite.

For this, we use the QR factorization $ Q^* \left[ \begin{array}{cc} D & G\end{array}\right] = \left[ \begin{array}{cc} \widehat{D} & \widehat{G}\\ 0_{k,  n} & 0_{k, r}\end{array} \right]$, where $\left[\begin{array}{cc}  \widehat{D} &  \widehat{G}\end{array}\right]$ has full row-rank and  $Q =\left[\begin{array}{cc}  Q_1&  Q_2 \end{array}\right]$ is unitary with $Q_2 \in \C^{n\times k}$. Set $D_{11} := Q_1^*DQ_1$ and  $ Q^*MQ =
\left[\begin{array}{cc} M_{11} & M_{12}  \\ M_{12}^* & M_{22}\end{array}\right]$, where $M_{22} \in \C^{k\times k}$.
{
Introduce  the matrix polynomials
\begin{eqnarray} \label{defH}H(\lam) &:=& \lam \left[\begin{array}{cc} M_{11} -  M_{12} M_{22}^{-1} M_{12}^* & 0 \\ 0  & -I_r\end{array}\right] + \left[ \begin{array}{cc} D_{11} &  \widehat{G} \\ {\widehat{G}}^* & 0 \end{array}\right], \\
\widehat{P}(\lam) &:=& \lam^2 (M_{11} -  M_{12} M_{22}^{-1} M_{12}^*) + \lam D_{11} +  \widehat{G}  \widehat{G}^*. \nonumber \end{eqnarray}}
{ Our next result shows that  the deflation procedure characterizes exactly how many nonzero eigenvalues remain after the deflation procedure.}
\begin{theorem}\label{new:thm2}  {  Consider  $ P(\lam) := \lam^2M+\lam D+GG^*$ and the associated  pencil $S(\lam)$ in \eqref{herm},} where $0<M=M^*\in \mathbb C^{n\times n}$,  $0\leq D=D^*\in \mathbb C^{n\times n}$, $ G \in \mathbb C^{{n\times r}}$ is of rank $r$.   Suppose that $\rank\left(\left[ \begin{array}{cc} D & G\end{array}\right]\right) = n-k$.

(a) Then  {  $ S(\lam)$ and $P(\lam)$ have exactly $ n+r -k$ nonzero eigenvalues counting multiplicity. In particular, $ 0 \in \Lam(S)$ and the algebraic multiplicity of $0$ is $k$. Also, $0 \in \Lam(P) $ and the algebraic multiplicity of $0$ is $n+k-r.$  }

(b) Let {  $H(\lam)$ and $\widehat P(\lam) $ be as in \eqref{defH}.}
Then, $\widehat{M} := M_{11} -  M_{12} M_{22}^{-1} M_{12}^* >0$ and, for all $ \lam \in \mathbb C$, there exists a nonsingular matrix $X \in \C^{(n+r) \times (n+r)}$ such that
\[ X^*S(\lam) X =  \diag
(H(\lam) , \lam M_{22})
\]
 and
 $ \Lam(\widehat{P})\setminus \{0\}= { \Lam(H)} = \Lam(S)\setminus\{0\} = \Lam(P)\setminus\{0\}$.


Additionally, if   $n= k+r$    then  $$ \Lam(\widehat{P})=\Lam(H) = \Lam(S)\setminus\{0\} = \Lam(P)\setminus\{0\}.$$

\end{theorem}

\proof Since { $ D^* = D$, we have $\rank\left(\left[ \begin{array}{cc} D & G\end{array}\right]\right) = \rank\left( \left[ \begin{array}{c} D \\ G^*\end{array}\right]\right) = n-k$.} Hence  $ \dim (\mathcal N(D)\cap \mathcal N(G^*)) = k$.  For $ x \in \mathcal N(D)\cap \mathcal  N(G^*)$, setting $ u := [ x^\top, 0 ]^\top \in \C^{n+k}$, we have $S(0) u = 0$ which shows that the multiplicity of $0 $ as an eigenvalue of $S(\lam)$ is at least $k$.

As $\rank(\left[\begin{array}{cc}  \widehat{D} &  \widehat{G}\end{array}\right]) = n-k$ and   $ G^*Q_2=0 = DQ_2$, it follows that the columns of $Q_2$ { belong to} $N(D)\cap N(G^*)$. Since $D^* = D$, we have $Q^*DQ =  \diag( D_{11}, 0)$.
Furthermore, we have
\[
\widehat{S}(\lam) :=\diag(Q^*, I_r) S(\lam) \diag(Q, I_r)= \lam \left[\begin{array}{cc|c} M_{11} & M_{12} & 0 \\ M_{12}^* & M_{22} & 0 \\ \hline  0 & 0 &  -I_r \end{array}\right] +  \left[\begin{array}{cc|c} D_{11} & 0 & \widehat{G} \\ 0 & 0 & 0 \\ \hline  {\widehat{G}}^* & 0 &  0 \end{array}\right].
\]
Now, consider the block  matrices
\[
E:= \left[\begin{array}{cc|c} I_{n-k} & 0 & 0 \\ 0 & 0 & I_k \\ \hline  0 & I_r &  0 \end{array}\right], \quad F :=  \left[\begin{array}{cc|c}I_{n-k}  & 0 & -M_{12}M_{22}^{-1}\\ 0  & I_r & 0 \\ \hline 0 & 0 & I_k\end{array}\right]^*.
\]
Then we have
\[
E^* \widehat{S}(\lam) E  =  \lam \left[\begin{array}{cc|c} M_{11}  & 0 & M_{12}\\ 0  & -I_r & 0 \\ \hline M_{12}^* & 0 & M_{22}\end{array}\right] + \left[ \begin{array}{cc|c} D_{11} &  \widehat{G} & 0 \\ {\widehat{G}}^* & 0 & 0 \\ \hline 0 & 0 & 0 \end{array}\right]\]
 and
 \begin{eqnarray*}
     F^*E^*\widehat{S}(\lam) EF &=& {\lam \left[\begin{array}{cc|c} M_{11}- M_{12}M_{22}^{-1}M_{12}^*  & 0 & 0\\ 0  & -I_r & 0 \\ \hline 0 & 0 & M_{22}\end{array}\right] + \left[ \begin{array}{cc|c} D_{11} &  \widehat{G} & 0 \\ {\widehat{G}}^* & 0 & 0 \\ \hline 0 & 0 & 0 \end{array}\right] } \\ &=& \diag( H(\lam), \lam M_{22}).
     \end{eqnarray*}
{ Hence, it follows that  $ \Lam(H) \setminus\{0\} = \Lam(S)\setminus\{0\} = \Lam(P)\setminus\{0\}$.}  By Lemma~\ref{lem:lem1} (see also Remark~\ref{rem:remtrim}), $H(\lam)$ is a trimmed linearization of $ \widehat{P}(\lam)$. Hence, $\Lam(H)\setminus\{0\} = \Lam(\widehat{P})\setminus\{0\}$.  Since $Q^*MQ$ is positive definite and $\widehat{M}$ is the Schur complement of $M_{22}$ in $Q^*MQ,$ it follows that $M_{11}- M_{12}M_{22}^{-1}M_{12}^*$ is positive definite.

Now we show that  $D_{11}$ is definite. Since  $\left[\begin{array}{cc}  \widehat{D} &  \widehat{G}\end{array}\right]$ has full row rank, we have $\mathcal N(D_{11})\cap \mathcal N({\widehat G}^*) =\{0\}$. Indeed, let $ x \in \mathcal N(D_{11})\cap \mathcal N({\widehat G}^*)$. Then $ D_{11}x =0$ and ${\widehat G}^* x =0$. Since $D$ is semidefinite, it follows from  $ x^*D_{11}x = 0 $ that $x^*Q_1 DQ_1x = 0$. This implies $DQ_1x = 0$ and hence $ x^* Q_1^*D =0 $ from which we obtain that $ x^*\widehat{D} =0$.
Since  $\left[\begin{array}{cc}  \widehat{D} &  \widehat{G}\end{array}\right]$ has full row rank, we get $ x =0$. Hence $D_{11}$ is definite, and therefore $\left[ \begin{array}{cc} D_{11} &  \widehat{G}  \\ {\widehat{G}}^* & 0  \end{array}\right]$ is nonsingular, and thus, $ 0\notin \Lam(H)$ and  we have $ \Lam(H) = \Lam(S)\setminus\{0\}$.

{ Since $\Lam(H) = \Lam(S)\setminus \{0\} = \Lam(P)\setminus\{0\}$ and  $\mathrm{Ind}(\mu, H) = \mathrm{Ind}(\mu, S) = \mathrm{Ind}(\mu, P) $ for all $\mu \in \Lam(H),$ it follows that $ H(\lam), S(\lam)$ and $P(\lam)$ have exactly the same number of nonzero eigenvalues counting multiplicity. As $H(\lam)$ has $n-k+r$ nonzero eigenvalues counting multiplicity, $S(\lam)$ and $P(\lam)$ have $n-k+r$ nonzero eigenvalues. Consequently, the algebraic multiplicity of $ 0 \in \Lam(S)$ is $ n+r - (n-k+r) = k$ and the algebraic multiplicity of $ 0 \in \Lam(P)$ is $ 2n -(n-k+r) = n+k-r.$}

Finally, note that $ \widehat{G} \in \C^{(n-k) \times r}$ and $\rank(\widehat{G}) = r$.  Since $ n= k+r$, it follows that $\widehat{G} \in \C^{r\times r}$ is nonsingular. Consequently, $0$ is not an eigenvalue of  $\widehat{P}(\lam) := \lam^2 (M_{11} -  M_{12} M_{22}^{-1} M_{12}^*) + \lam D_{11} +  \widehat{G}  \widehat{G}^*$.  Hence, $ \Lam(\widehat{P}) = \Lambda(H) = \Lambda(P)\setminus\{0\}.$
\eproof

{
\begin{remark}\label{rem:schurcompl}{\rm
Using the Schur complement in the deflation procedure may lead to the question as to what happens if the matrix $M_{22}$ is ill-conditioned, and whether this has a negative effect on the computation of eigenvalues in the reduced problem. This question has been addressed in a slightly different context in the deflation of infinite eigenvalues in symmetric/skew-symmetric pencils \cite{MehX14a}, where by a detailed perturbation analysis it is shown that using the Schur complement approach in structure preserving deflation methods gives results that are as accurate as using methods based on unitary  projections. Below, in Example~\ref{ex:Schur}  we present a numerical  test to illustrate this observation. }
\end{remark}
}

\begin{remark} \label{rem:trimmed2}{\rm { Since for the matrix polynomials in \eqref{defH} $H(\lam)$ is a trimmed linearization of $\widehat{P}(\lam)$ { with  $\Lam(H) = \Lam(\widehat{P})\setminus\{0\}$,  we have  $\mathrm{Ind}(\mu, H) = \mathrm{Ind}(\mu, \widehat{P})$ for all $\mu \in  \Lam(H)$.} It follows from Theorem~\ref{new:thm2} that $\mathrm{Ind}(\mu, H) = \mathrm{Ind}(\mu, S) $ for all $\mu \in { \Lam(H)}$. Consequently,
we have  $\mathrm{Ind}(\mu, \widehat{P})= \mathrm{Ind}(\mu, H) = \mathrm{Ind}(\mu, S) = \mathrm{Ind}(\mu, P) $ for all $\mu \in \Lam(P)\setminus\{0\}$. Hence $H(\lam)$ is also a trimmed linearization of $P(\lam)$    and   $H(\lam)$ deflates the zero eigenvalues of $P(\lam)$.} }
\end{remark}

\begin{remark} \label{rem:reduced}{\rm
Theorem~\ref{new:thm2} shows that we can always construct a trimmed linearization $H(\lam)$  that deflates the zero eigenvalues of $P(\lam)$ when the damping matrix is semidefinite. Additionally, if  $n= k+r$  then $H(\lam)$ yields a structure-preserving reduced size quadratic polynomial $\widehat{P}(\lam)$ in which the zero eigenvalues of $P(\lam)$ have been completely deflated.}
\end{remark}

 We now show 
 that the assumption that $D$ is semidefinite cannot be omitted.

\begin{example} \label{ex:nodrop}{\rm Consider $ P(\lam) := \lam^2 I + \lam D + GG^*,$ where $G := \left[\begin{array}{cc} 1 & 0 \end{array}\right]^\top $ and $D := \left[\begin{array}{cc} 1 & 1 \\ 1 & 0\end{array}\right]$ is indefinite. Then
$
S(\lam) =\lam \left[\begin{array}{cc|c} 1 & 0 & 0 \\ 0 & 1 & 0 \\ \hline  0 & 0 &  -1 \end{array}\right] +  \left[\begin{array}{cc|c} 1 & 1 & 1 \\ 1 & 0 & 0 \\ \hline  1 & 0 &  0 \end{array}\right]
$.
Note that  $\rank(\left[\begin{array}{cc}  D &  G\end{array}\right]) = 2$ and clearly $ 0 \in \Lam(S)$. This shows that Lemma~\ref{lem:new2} may not hold if $D$ is not semidefinite.

Also note that $G^*D^{-1}G = 0$
and for $ Q_2 := e_2$ we have $Q_2^* G = 0$ and
  $ Q_2^*DQ_2 = 0$. Thus the conditions in the other results are not satisfied either.
  }
\end{example}

\subsection{Computational methods} We now briefly discuss two  methods for the implementation of the result in Theorem~\ref{new:thm2}. {\sc Matlab} implementations are presented in \cite{AlaMT25_arxiv}. Consider the matrices $ 0 < M \in \C^{n\times  n}$, $D = D^* \in \C^{n\times  n}$, and $ G \in \C^{n\times  r}$. Suppose that $D$ is semidefinite and that $\rank(G) = r$. Also, suppose that  $m : = \rank( \left[\begin{array}{cc} D & G \end{array}\right]) < n. $ \\

For the first method, consider a rank revealing QR factorization
\begin{equation}\label{rrqr}
\left[\begin{array}{cc} D & G \end{array}\right] \Pi = Q  \left[\begin{array}{cc} R_{11} & R_{12} \\ 0 & 0 \end{array}\right],
\end{equation}
where $ R_{11} \in \C^{m \times m}$ is nonsingular, $Q \in \C^{n\times n}$ is unitary and $\Pi \in \C^{(n+r)\times  (n+r)}$ is a permutation matrix.   Partition $ Q =\left[\begin{array}{cc} Q_1 & Q_2 \end{array}\right]$, where $ Q_1 \in \C^{n\times  m}.$ Then it follows that $ Q_2^* \left[\begin{array}{cc} D & G \end{array}\right] = 0$ which shows that   $\mathrm{span}(Q_2) = \mathcal{N}(D) \cap \mathcal{N}(G^*)$.  Also,    $Q_1^*\left[\begin{array}{cc} D & G \end{array}\right] = \left[\begin{array}{cc} R_{11} & R_{12} \end{array}\right]\Pi^* $ has full row rank. Hence we have the following algorithm. \\

\begin{algorithm} \label{alg1} QR-based method for deflating zero eigenvalues of $P(\lambda)$.
 \begin{itemize}
 	\item[1.] Compute $m : = \rank( \left[\begin{array}{cc} D & G \end{array}\right]).$ If $ m = n$ then STOP as the zero eigenvalue is already deflated. If $ m < n$ then proceed as follows.
 	\item[2.]  Compute a rank revealing QR factorization as in \eqref{rrqr}.
 	\item[3.] Partition $ Q =\left[\begin{array}{cc} Q_1 & Q_2 \end{array}\right]$  conformably and compute $ \widehat{D} := Q_1^*DQ_1,$ where $ Q_1 \in \C^{n\times  m}$.
 	\item[4.]  Compute $ \left[\begin{array}{cc} M_{11} & M_{12} \\ M_{12}^* & M_{22} \end{array}\right]  := Q^*MQ$ and the Schur complement $ \widehat{ M} := M_{11}- M_{12}M_{22}^{-1} M_{12}^*$, where $ M_{11} \in \C^{m \times m}$.
 	\item[5.]  Compute $ R :=  \left[\begin{array}{cc} R_{11} & R_{12} \end{array}\right]\Pi^*$ and define $\widehat{G}$ to be the last $r$ columns of $R$, where $r = \rank(G)$.
 	\item[6.]  Construct the pencil $ \lam \left[\begin{array}{cc} \widehat{M} & 0 \\ 0 & -I_r \end{array}\right] + \left[\begin{array}{cc} \widehat{D} & \widehat{G} \\ \widehat{G}^* & 0 \end{array}\right]$ in which  the zero eigenvalues are deflated.
 \end{itemize}
\end{algorithm}

The second method is based on the singular value decomposition (SVD)
\begin{equation}\label{svddg}
\left[\begin{array}{cc} D & G \end{array}\right]  = U \diag(\Sigma_m, 0)V^*,
\end{equation}
where $ \Sigma_m \in \C^{m\times m}$ is diagonal, containing all $m$ nonzero singular values on the diagonal, and  $ U \in \C^{n\times  n}$, $V \in \C^{(n+r) \times (n+r)}$  are unitary. Partition $ U =\left[\begin{array}{cc} U_1 & U_2 \end{array}\right]$, where $ U_1 \in \C^{n\times  m}$. Then it follows that $ U_2^* \left[\begin{array}{cc} D & G \end{array}\right] = 0$ which shows that $\mathrm{span}(U_2) = \mathcal{N}(D) \cap \mathcal{N}(G^*)$.  Also,    $U_1^*\left[\begin{array}{cc} D & G \end{array}\right] = \left[\begin{array}{cc} \Sigma_m & 0 \end{array}\right]V^* $ has full row rank. This yields the following algorithm.
	
\begin{algorithm}\label{alg2} SVD-based method for deflating zero eigenvalues of $P(\lambda)$.
		\begin{itemize}
				\item[1.] Compute $m : = \rank( \left[\begin{array}{cc} D & G \end{array}\right]).$ If $ m = n$ then STOP as the eigenvalue $0$ is already deflated. If $ m < n$ then proceed as follows.
				\item[2.]  Compute the SVD as in \eqref{svddg}.
		\item[3.] Partition $ U =\left[\begin{array}{cc} U_1 & U_2 \end{array}\right]$  and compute $ \widehat{D} := U_1^*DU_1,$ where $ U_1 \in \C^{n\times  m}.$
		\item[4.]  Compute $ \left[\begin{array}{cc} M_{11} & M_{12} \\ M_{12}^* & M_{22} \end{array}\right]  := U^*MU$ and $ \widehat{ M} := M_{11}- M_{12}M_{22}^{-1} M_{12}^*\in \C^{m \times m}$.
		\item[5.]  Partition $ V =\left[\begin{array}{cc} V_1 & V_2 \end{array}\right]$ and compute $ R :=  \left[\begin{array}{cc} \Sigma_m &  0 \end{array}\right]V^* = \Sigma_m V_1^*$, where $ V_1 \in \C^{(n+r) \times m}$.
		 \item[6.] Define $\widehat{G}$ to be last $r$ columns of $R,$ where $r = \rank(G)$.
		\item[7.]  Construct the pencil $ \lam \left[\begin{array}{cc} \widehat{M} & 0 \\ 0 & -I_r \end{array}\right] + \left[\begin{array}{cc} \widehat{D} & \widehat{G} \\ \widehat{G}^* & 0 \end{array}\right]$ in which the $0$ eigenvalues have been deflated.
	\end{itemize}
	\end{algorithm}
    \begin{remark}
        \label{both0andinf}{\rm  It has been shown in \cite{MacMMM06b} that there is no structured linearization if both $M$ and $K$ are singular. In this case the discussed methods cannot be employed directly without first projecting out one of the two eigenvalues via an unstructured linearization.}
    \end{remark}

{
\subsection{Numerical Examples}\label{sec:example}
In this subsection we illustrate the two  methods described in the previous subsection via several numerical examples.
The {\sc Matlab} codes that have been used are presented in the Appendix of \cite{AlaMT25_arxiv}.
Our first example demonstrates the procedures for a large scale problem.
\begin{example}\label{ex:nino}{ \rm
Consider a problem arising in the free vibration of a rod  as presented in \cite{TambacaTruhar2025}. A detailed description of the problem is given in \cite{AlaMT25_arxiv}. The resulting
matrices $M$ and $K$ have size $ n = 1602$, with
singular $M$.

Applying our methods, we deflate the infinite eigenvalues by deflating the zero eigenvalues of the reverse polynomial
\[
\operatorname{rev}P(\mu) =
\mu^2 K + \mu D + M \equiv \mu^2 \tilde M + \mu D + \tilde G \tilde G^*,
\]
with $ M = \tilde G \tilde G^* $, $ \tilde G \in \mathbb{R}^{n \times r} $, $ r = 1599 $ and $\tilde M$ nonsingular.

Applying the {\sc Matlab} function  \texttt{eig}  to the pair $ (\tilde M, \tilde G\tilde G^*) $, yields three eigenvalues $ 0 $, $ 8.351 \times 10^{-22} $, and $ 1.0829 \times 10^{-19} $  that are numerically zero and that will  have to be deflated.

The damping matrix is defined as:
\[
D = \sum_{i \in \mathcal{I}_d} v_i d_i d_i^T,
\]
where $ d_i $ are the selected rows of the eigenvector matrix of the pair $ (\tilde M, \tilde G\tilde G^*) $ corresponding to the damper locations in $ \mathcal I_d = \{124, 125, 126, 127, 144\} $, with $ r_d = | \mathcal I_d | = 5 $, and we set all viscosities to $ v_i = 10 $.

A {\sc Matlab} calculation gives
$
\operatorname{rank}(\tilde G) = 1599$, $ m = \operatorname{rank}([D, \tilde G]) = 1601$
and thus $k=1$.
In the calculated pencil
$$
\lambda \begin{bmatrix} \widehat{M} & 0 \\ 0 & -I_r \end{bmatrix} + \begin{bmatrix} \widehat{D} & \widehat{G} \\ \widehat{G}^* & 0 \end{bmatrix}
$$
of size $ 3200 \times 3200 $ the zero eigenvalues are deflated,  the smallest eigenvalue modulus is $|\lambda_{\min}| = 1.211$, and  $\max_i({Re}(\lambda_i)) = -0.4598$.
}
\end{example}

Our next example illustrates Theorem~\ref{new:thm2}.
%
\begin{example}\label{ex:exactevs}{\rm
In this example, the matrices $(M, D, G)$  are generated randomly  via the {\sc Matlab} code presented in \cite{AlaMT25_arxiv} with input {\tt n= 5, r =3,  d = [2, 0, 2, 0, 0]'; w = [ 2, 3, 4 ]'; }. In this case, $ \rank(D) =  2$, $\rank(G) = 3$ and $m:=\rank[( D, G]) = 5$.  Hence $P(\lambda)$ has $m+r =8$ nonzero eigenvalues and $2n-8 = 2$ zero eigenvalues.  	The computed eigenvalues of $P(\lambda)$ by the {\sc Matlab} code {\tt  polyeig(K, D, M)} are given below.
	\begin{verbatim}	
		-7.263133603723262e-03 + 6.081955067660005e-01i
		-7.263133603723262e-03 - 6.081955067660005e-01i
		-1.248073295181371e-02 + 5.508782232131781e-01i
		-1.248073295181371e-02 - 5.508782232131781e-01i
		-5.518062238563429e-02 + 3.385503421073031e-01i
		-5.518062238563429e-02 - 3.385503421073031e-01i
		-2.752279210033463e-01 + 0.000000000000000e+00i
		-1.238978644744203e-03 + 0.000000000000000e+00i
		 1.507566034107082e-15 + 1.627212285681345e-15i
		 1.507566034107082e-15 - 1.627212285681345e-15i
		\end{verbatim}
Note that the last two eigenvalues are zero eigenvalues and the pencil $S(\lambda)$ has $n-m =0$ zero eigenvalues. Hence $S(\lambda)$ deflates all the zero eigenvalues of $P(\lambda)$. The eigenvalues of $S(\lambda)$ computed by {\tt eig} are as follows.
	
	\begin{verbatim}		
		-7.263133603722457e-03 + 6.081955067660011e-01i
		-7.263133603722457e-03 - 6.081955067660011e-01i
		-1.248073295181332e-02 + 5.508782232131768e-01i
		-1.248073295181332e-02 - 5.508782232131768e-01i
		-5.518062238563402e-02 + 3.385503421073032e-01i
		-5.518062238563402e-02 - 3.385503421073032e-01i
		-1.238978644741443e-03 + 0.000000000000000e+00i
		-2.752279210033460e-01 + 0.000000000000000e+00i
	\end{verbatim}

}
	
\end{example}

The next example illustrates that the calculation of the Schur complement is not critical for the eigenvalues of the deflated problem.
\begin{example}\label{ex:Schur}    {\rm
In this example, the matrices $(M, D, G)$  are generated randomly  via the {\sc Matlab} code presented in \cite{AlaMT25_arxiv} with input  {\tt n= 5, r =2,  d = [2, 0, 2, 0, 0]'; w = [ 2, 3 ]'; } In this case, $ \rank(D) = \rank(G) =2$ and $m:=\rank[( D, G]) = 4.$  Hence $P(\lam)$ has $m+r =6$ nonzero eigenvalues and $2n-6 = 4$ zero eigenvalues. The eigenvalues of $P(\lam)$ computed by {\tt  polyeig(K, D, M)} are given below.
\begin{verbatim}	
	-8.354874652713031e-03 + 6.092645960144358e-01i
	-8.354874652713031e-03 - 6.092645960144358e-01i
	-5.537232363234557e-02 + 3.376550957779614e-01i
	-5.537232363234557e-02 - 3.376550957779614e-01i
	-2.929077221606973e-01 + 0.000000000000000e+00i
	-5.953758799614365e-03 + 0.000000000000000e+00i
	-2.245340872846755e-16 + 0.000000000000000e+00i
	 2.654445177086020e-09 + 0.000000000000000e+00i
	-2.654444863332211e-09 + 0.000000000000000e+00i
	 9.536161154148316e-16 + 0.000000000000000e+00i
	\end{verbatim}
The last four eigenvalues arise from zero eigenvalues and the pencil $S(\lam)$ has  $n-m = 1$ zero eigenvalue. The eigenvalues of $S(\lam)$ computed by {\tt eig} are given below. Note that last but one eigenvalue is the zero eigenvalue of $S(\lam)$.
\begin{verbatim}
-8.354874652713066e-03 + 6.092645960144352e-01i
-8.354874652713066e-03 - 6.092645960144352e-01i
-5.537232363234516e-02 + 3.376550957779615e-01i
-5.537232363234516e-02 - 3.376550957779615e-01i
-2.929077221606972e-01 + 0.000000000000000e+00i
-3.398002201242367e-18 + 0.000000000000000e+00i
-5.953758799613264e-03 + 0.000000000000000e+00i
\end{verbatim}
%

The  pencil $H(\lam)$ deflates the zero eigenvalue of $S(\lam)$. The eigenvalues of $H(\lam)$ computed by {\tt eig} are given below.
%
\begin{verbatim}		
	-8.354874652713024e-03 + 6.092645960144357e-01i
	-8.354874652713024e-03 - 6.092645960144357e-01i
	-5.537232363234555e-02 + 3.376550957779611e-01i
	-5.537232363234555e-02 - 3.376550957779611e-01i
	-5.953758799613271e-03 + 0.000000000000000e+00i
	-2.929077221606972e-01 + 0.000000000000000e+00i
	\end{verbatim}
%
The QR-based method produces exactly similar results.
}
\end{example}

}
\section{{ Hyperbolic eigenvalue problems}} \label{sec:hyperbolic}
The solution of quadratic hyperbolic eigenvalue problems  of the form \eqref{qep} satisfying that $M>0$ and \eqref{defhyp} is a widely investigated problem, see, e.g., \cite{Ple06,Sap76,Ves11}. In this section we study trimmed linearizations  of hyperbolic  quadratic eigenvalue problems and also revisit the inertia. For this, we make use of an equivalent $J$-Hermitian (also called pseudo-Hermitian)  standard eigenvalue problem which is obtained from the  pencil $L(\lam)$ given in \eqref{jherm}.

Let $ M = CC^*$ be the Cholesky decomposition of $M$.  Define
\begin{equation}\label{eq:std} \mathcal{A} := \left[\begin{matrix}
	-C^{-1}D (C^{-1})^* & -C^{-1} G\\  (C^{-1} G)^* & 0
\end{matrix}\right], \  J := \left[\begin{matrix}
-I_n & \\ & I_p
\end{matrix}\right].  \end{equation}
Then $\cA$ is \emph{$J$-Hermitian}, i.e., $(J\cA)^* = J\cA$,  and since $ JL(\lam) = -S(\lam)$,  we have
\begin{equation}\label{eq:jpl} L(\lam) = \left[\begin{matrix} C & \\ & I_r
\end{matrix} \right] (\lam I_{n+r} - \mathcal{A}) \left[ \begin{matrix} C^* & \\ &  I_r
\end{matrix}\right],
S(\lam) = \left[\begin{matrix} C & \\ & I_r
\end{matrix} \right]  (J\mathcal{A} - \lam J) \left[ \begin{matrix} C^* & \\ &  I_r
\end{matrix}\right].
\end{equation}
We refer to $\cA$ as the \emph{standard $J$-Hermitian trimmed linearization} of $P(\lam)$.
Let $\Lam (\cA)$ denote the spectrum of $\cA$. Then  by Lemma~\ref{lem:lem1}, we have   $\Lambda(\cA) \setminus \{0\} = \Lambda(P)\setminus\{0\}$ and, by Lemma~\ref{lem:new1}, we have  $\Lambda(\cA) = \Lambda(P)\setminus\{0\}$, whenever $D$ is semidefinite and  $\rank(\left[\begin{array}{cc}  D &  G\end{array}\right]) = n$. \\

We mention that a permuted version of $\mathcal{A}$,   referred to as the phase-space matrix, has been obtained in~\cite[p.24]{Ves11} by employing an  entirely different method in the case when both $M$ and $K$ are positive definite.

\begin{definition}\cite{Ves11}\label{def:jherm}
	A nonzero $ u \in \C^{n+r}$  is said to be \emph{$J$-neutral} if $ u^*Ju =0$, and \emph{$J$-positive} (resp., \emph{$J$-negative}) if  $ u^*Ju  >0$ (resp.,  $ u^*Ju <0$).  A subspace $ V \subset \C^{n+r}$ is said to be \emph{$J$-neutral} (resp., \emph{$J$-positive , $J$-negative}) if each nonzero $ u \in V$ is $J$-neutral (resp., $J$-positive, $J$-negative).
\end{definition}

Since $\cA$ is $J$-Hermitian, the eigenvalues of $\cA$ are distributed in the complex plane symmetrically with respect to the real axis. Furthermore, if $v$ is an eigenvector of $\cA$ corresponding to a complex  eigenvalue of $\cA$ with nonzero imaginary part then  $v$ is $J$-neutral. Indeed, let $ \cA v= \lam  v$ and let $ \lam$ be complex with nonzero imaginary part. Then $ v^* J\cA v = \lam v^*Jv$ shows that $v^*Jv = 0 $ as both $J$ and $ J\cA$ are Hermitian.

Let $ \lam $ be a real eigenvalue of $\cA$. Then $ \lam$ is  $J$-positive  (resp., $J$-negative) if the spectral subspace $ \mathcal{E}(\lam) :=\mathcal N\left((\cA - \lam I_{n+r})^{n+r}\right)$ is $J$-positive (resp., $J$-negative).
(Note that sometimes  the terminology \emph{positive type (resp., negative type)} is used.)
If the spectral subspace $ \mathcal{E}(\lam)$ contains a nonzero $J$-neutral vector then $ \lam$ is said to be of \emph{mixed type (or $\lam$ has mixed sign characteristics)}. Equivalently, let $ X \in \C^{(n+r) \times m}$ be a full column rank matrix such that $\mathrm{span}(X) = \mathcal{E}(\lam)$.
Then  $ \lam$  is $J$-positive (resp., $J$-negative) if and only if $ X^*JX >0$ (resp., $X^*JX < 0$). Finally,  $\lam$ is of mixed type if and only if $ X^*JX$ is indefinite.
We say that $\cA$ has {\em definite spectrum} (or that the spectrum of $\cA$ is definite) if each eigenvalue of $\cA$ is either $J$-positive or $J$-negative. \\

A $J$-Hermitian matrix  with  definite spectrum  admits a special decomposition and the eigenvalues have special properties; see~\cite{Ves11}. For reference, we summarize some of these results  for $\cA$ in the following theorem.   \\

\begin{theorem} \label{cf} Let $ \lam $ be a $J$-positive or  $J$-negative eigenvalue of $\cA$ in \eqref{eq:std} of algebraic multiplicity $m$. Then there exists a nonsingular matrix $U$ such that $$ U^{-1} \cA U = \diag( \lam I_m,\widehat{A})\; \text{ and } \; U^* JU =  \diag(\epsilon I_m , \widehat{J}  ),$$ where $ (\widehat{J}\widehat{A})^*= \widehat{J}\widehat{A}$ and $ \lam \notin \Lam(\widehat{A})$.  Here $ \epsilon =1$ if $\lam$ is $J$-positive and $\epsilon = -1$ if $\lam$ is $J$-negative. In particular, if $\cA$ has definite spectrum with distinct eigenvalues $\lam_1, \ldots, \lam_\ell$ and multiplicities $m_1, \ldots, m_\ell$, respectively, then there exists a nonsingular matrix $U$ such that
	$$ U^{-1}\cA U = \lam_1 I_{m_1} \oplus \cdots \oplus \lam_\ell I_{m_\ell} \text{ and } U^*JU =  \epsilon_1 I_{m_1} \oplus \cdots \oplus \epsilon_\ell I_{m_\ell},$$ where $ \epsilon_j = 1$ if $ \lam_j$ is $J$-positive  and $\epsilon_j=-1$ if $\lam_j$ is $J$-negative for $ j=1,\ldots,\ell$. Furthermore,  the spectrum of $\cA$ consists of $ n$ (counting multiplicity) $J$-negative and $r$ (counting multiplicity) $J$-positive eigenvalues.
\end{theorem}


We now show that the Hermitian pencil $S(\lam)$ in (\ref{herm}) is  a definite trimmed linearization of $P(\lam)$ when $P(\lam)$ is overdamped, i.e., if $D>0$.

\begin{theorem}\label{th:dp} { Consider  $ P(\lam) := \lam^2M+\lam D+GG^*$ and the associated pencil $S(\lam)$ in \eqref{herm},}  where $0<M=CC^*\in \mathbb C^{n\times n}$,  $0<D=D^*\in \mathbb C^{n\times n}$, $ G \in \mathbb C^{{n\times r}}$ is of rank $r$.
Then $ P(\lam)$ is hyperbolic if and only if $ S(\lam)$ is a definite pencil.
For $ \mu \in \R\setminus \{0\}$,  $ \imath(S(\mu)) = \imath(P(\mu)) + ( 0, r, 0)$ when $ \mu >0$ and $ \imath(S(\mu)) = \imath(-P(\mu)) + ( r, 0, 0)$ when $ \mu < 0$.
\end{theorem}

\proof Suppose that $ P(\lam)$ is hyperbolic. Then  $-P(\mu) >0 $ for some real $ \mu$ and thus $ \mu <0$, since if  $\mu= 0$ then $ P(\mu) = GG^* \geq 0$ contradicting that $P(\mu) <0$. Similarly, if $ \mu >0$  then $P(\mu) >0$ which contradicts that $P(\mu) <0$.

Now consider the matrix $ S(\mu)$.  Then $ { T(\mu)} := P(\mu)/\mu =\mu M + D + GG^*/\mu $ is the Schur complement of $ -\mu I_r$ in $S(\mu)$ and
\begin{equation}\label{schurcomp} S(\mu) = \left[\begin{matrix}
I_n &  -G/\mu \\ & I_r
\end{matrix} \right] \left[\begin{matrix}
  T(\mu) &   \\ &  -\mu I_r
\end{matrix} \right] \left[\begin{matrix}
I_n &  -G/\mu \\ & I_r
\end{matrix} \right]^*.
\end{equation}

Thus, $ S(\mu) > 0$ if and only if $ -\mu I_p >0$ and $T(\mu) >0$. Since $ \mu <0,$ we have $ -\mu I_p >0$ and $ T(\mu) = P(\mu)/\mu = \frac{- P(\mu)}{-\mu} >0$. Consequently, $S(\mu)$ is  positive definite and hence the pencil $S(\lam)$ is definite. Conversely, suppose that $ S(\mu) > 0$ for some real $ \mu$. Then $- \mu I_r >0 $ and $ T(\mu) >0$.   Hence we have $ \mu <0$ and $ - P(\mu) = -\mu  T(\mu) >0$.  This shows that $ P(\lam)$ is hyperbolic.

By (\ref{schurcomp}), we have $ \imath(S(\mu)) = \imath(T(\mu)) + \imath(-\mu I_p)$. Then the   assertions about the inertia follows from   the fact that $ \imath(T(\mu)) = \imath(P(\mu))$ when $ \mu > 0$ and $ \imath(T(\mu)) = \imath(-P(\mu))$ when $ \mu < 0$.
\eproof


{ \begin{remark}  Note that $S(\lam)$ in Theorem~\ref{th:dp} is a  definite trimmed linearization of $P(\lam)$  with $ \Lambda(S) = \Lambda(P)\setminus\{0\}$ and that the assumption $ D >0$ cannot be dropped.  If $D$ is positive semidefinite  and $ N(G^*) \subsetneq N(D)$ then $P(\lam)$ in Theorem~\ref{th:dp} is not  hyperbolic as $(x^*Dx)^2 - 4 x^*xMx x^* GG^*x  < 0$ for nonzero $ x \in N(D)\setminus N(G^*).$ On the other hand, if $ N(D) \subset N(G^*)$ then $ P(\lambda)$ can be weakly hyperbolic, i.e., $(x^*Dx)^2 - 4 x^*xMx x^* GG^*x \geq 0$ for nonzero  $ x$.  We mention that definite linearizations of a hyperbolic $P(\lam)$  exist without the requirements that $D\geq 0$ and $K\geq 0$, see~\cite{HigMT09,Ple06}.
\end{remark}  }

A hyperbolic quadratic matrix polynomial has many interesting properties, see~\cite{HigMT09,Ple06}, such as
\begin{itemize}
\item [(a)] all its eigenvalues are  real and semisimple,
\item [(b)] there is a gap between the largest $n$ and  the smallest $n$ eigenvalues, and
\item [(c)]
the inertia of the hyperbolic matrix polynomial  at $\mu \in \R$ yields the number of eigenvalues larger/smaller than $\mu$.
\end{itemize}

We now show similar results for the $J$-Hermitian matrix $\cA$.
Let $ \Lambda_1$ and $ \Lambda_2$ be finite subsets of $\R$ and let $ \alpha \in \R$.  We write $ \Lambda_1 < \alpha < \Lambda_2$ if $ \max\{ \lam : \lam \in \Lambda_1\}  < \alpha < \min\{\mu : \mu \in \Lambda_2\}$.  Then we have the following result.

\begin{theorem} \label{th:defspec} Let $ P(\lam) := \lam^2M+\lam D+GG^*$ be  hyperbolic, where $0<M=CC^*\in \mathbb C^{n\times n}$,  $0<D=D^*\in \mathbb C^{n\times n}$, $ G \in \mathbb C^{{n\times r}}$ is of rank $r$.
Then we have the following assertions.
\begin{itemize}
\item[(a)] The spectrum of $\cA$ has $ n$ (counting multiplicity) $J$-negative and $r$ (counting multiplicity) $J$-positive eigenvalues.

\item[(b)] Let $ \Lambda^{+}(\cA)$ and $\Lambda^{-}(\cA),$ respectively, denote the set of $J$-positive and $J$-negative eigenvalues of $\cA$. Then there exists $ \alpha \in \R$ such that $ \Lambda^{-}(\cA)< \alpha < \Lambda^+(\cA)$.

\item[(c)] Let $ \lam^+_{\min}$ and $ \lam^{-}_{\max},$ respectively, denote the smallest and the largest eigenvalues in $ \Lambda^{+}(\cA)$ and $\Lambda^{-}(\cA)$. Then $ (\lam^{-}_{\max},\,\, \lam^+_{\min})$ is the definitizing interval for $\cA$ in the sense that for each $ \mu \in  (\lam^{-}_{\max},\,\, \lam^+_{\min})$ the matrix $ J\cA - \mu J$ is positive definite.

\item[(d)] For each $ \mu \in  (\lam^{-}_{\max},\,\, \lam^+_{\min})$, we have $ \mu< 0$, $S(\mu) >0$ and $-P(\mu) >0$, where $S(\lam)$ is as in (\ref{herm}).

\item[(e)] Suppose  that $ \mu \in \R $ is not an  eigenvalue of $\cA$. If $\mu > \lam^+_{\min}$ then  we have
\begin{eqnarray*} \imath_{+}(J\cA- \mu J) &=&  \#\{ \lam \in \Lambda^+(\cA) : \lam > \mu\}+n , \\ \imath_{-}(J\cA -\mu J) &=& \#\{ \lam \in \Lambda^+(\cA) : \lam < \mu\}.
\end{eqnarray*}
On the other hand, if $ \mu < \lam^{-}_{\max}$ then we have
\begin{eqnarray*} \imath_{+}(J\cA- \mu J) &=&  \#\{ \lam \in \Lambda^{-}(\cA) : \lam < \mu\}+r , \\ \imath_{-}(J\cA -\mu J) &=& \#\{ \lam \in \Lambda^{-}(\cA) : \lam > \mu\}.
\end{eqnarray*}
\end{itemize}
\end{theorem}

\proof By (\ref{eq:jpl}), we have
\[
S(\lam) = \diag( C, I_r)\
 ( J\mathcal{A} -\lam J )\ \diag( C^*,  I_r),
\]
and by Theorem~\ref{th:dp},  we have that $S(\lam)$ is a definite pencil. Hence there exists $\mu \in \R$ such that $ S(\mu) >0 $ which implies that $ J\cA -J\mu  >0$. This shows that  $\phi(z) :=  z- \mu$ is a definitizing polynomial of $\cA$, i.e., $ J\phi(A) = J\mu- J\cA >0$. By \cite[Theorem~10.3]{Ves11}, then $\cA$ has definite spectrum. Hence, by Theorem~\ref{cf}, $\Lam(\cA)$ has $r$ $J$-positive eigenvalues and $n$ $J$-negative eigenvalues.

As $\cA$ is definitizable, i.e.,  $JA -J \mu >0$,  the assertion (b)  follows from \cite[Theorem~10.6]{Ves11}.  By Theorem~\ref{cf}, we have
\[
U^{-1}\cA U = \lam_1 I_{m_1} \oplus \cdots \oplus \lam_\ell I_{m_\ell} \text{ and } U^*JU =  \epsilon_1 I_{m_1} \oplus \cdots \oplus \epsilon_\ell I_{m_\ell},
\]
where $ \lam_1, \ldots, \lam_\ell$ are the distinct eigenvalues of $\cA$ with multiplicities $m_1, \ldots, m_\ell$, respectively, and $\epsilon_j \in \{1, -1\}$ for $ j=1,\ldots,\ell$. Let $ \mu \in  (\lam^{-}_{\max},\,\, \lam^+_{\min})$. Then  $\epsilon_j(\lam_j - \mu) >0 $ for $j=1,\ldots,\ell$. Hence
\begin{eqnarray}\label{eigsign} U^*(J\cA- J \mu )U &=&  U^*JU U^{-1} \cA U  -  U^*JU \mu\\ \nonumber &=& \epsilon_1(\lam_1 -\mu)  I_{m_1} \oplus \cdots \oplus \epsilon_\ell(\lam_\ell -\mu)I_{m_\ell} > 0 \end{eqnarray}
which proves (d).

If $ \mu \in  (\lam^{-}_{\max},\,\, \lam^+_{\min})$ then   $ J\cA - \mu J > 0$. Hence by (\ref{eq:jpl}), we have that $S(\mu) >0. $ Consequently,  by (\ref{schurcomp}) we have $ \mu <0$ and $-P(\mu) >0$. This proves (c).

Finally, (e) follows from (\ref{eigsign}). Indeed, if $ \mu \in \R$ is not an eigenvalue of $\cA$ then
\begin{equation}\label{inert} \imath(J\cA - \mu J) = \imath\left( \epsilon_1(\lam_1 -\mu)  I_{m_1} \oplus \cdots \oplus \epsilon_\ell(\lam_\ell -\mu)I_{m_\ell}  \right).\end{equation}
Thus, if $\mu > \lam^+_{\min}$ then  $\epsilon_j(\lam_j- \mu) >0$ for each $J$-negative  eigenvalue $\lam_j$. As  $\cA$ has $n$ $J$-negative eigenvalues, it follows  from (\ref{inert})
 that
 \begin{eqnarray*} \imath_{+}(J\cA- \mu J) &=&  \#\{ \lam \in \Lambda^+(\cA) : \lam > \mu\}+n , \\ \imath_{-}(J\cA -\mu J) &=& \#\{ \lam \in \Lambda^+(\cA) : \lam < \mu\}.
 \end{eqnarray*}
 Similarly,  if $ \mu < \lam^{-}_{\max}$   then  $\epsilon_j(\lam_j- \mu) >0$ for each $J$-positive  eigenvalue $\lam_j$. As $\cA$ has  $r$ $J$-positive eigenvalues, it follows from (\ref{inert}) that
 \begin{eqnarray*} \imath_{+}(J\cA- \mu J) &=&  \#\{ \lam \in \Lambda^{-}(\cA) : \lam < \mu\}+r , \\ \imath_{-}(J\cA -\mu J) &=& \#\{ \lam \in \Lambda^{-}(\cA) : \lam > \mu\}.
 \end{eqnarray*}
This completes the proof.
\eproof

\begin{remark} \label{rem:congru} {\rm

(a) Let $P(\lam) := \lam^2 M+ \lam D+ GG^*$ be Hermitian with $ M=CC^*>0$ and $ \rank(G) = r$. 	By (\ref{eq:jpl}) and (\ref{schurcomp}), we have that
	 \begin{eqnarray}\label{congru}
     S(\lam)&=&\left[\begin{matrix} C & \\ & I_r
	 \end{matrix} \right]  (J\mathcal{A} - \lam J) \left[ \begin{matrix} C^* & \\ &  I_r
	 \end{matrix}\right]\\ \nonumber  & =& \left[\begin{matrix}
		I_n &  -G/\lam \\ & I_r
	\end{matrix} \right] \left[\begin{matrix}
		P(\lam)/\lam &   \\ &  -\lam I_r
	\end{matrix} \right] \left[\begin{matrix}
		I_n &  -G/\lam \\ & I_r
	\end{matrix} \right]^*.
    \end{eqnarray}	
    It follows that $ J\cA - \mu J >0$ for some $ \mu \in \R$ if and only if $\mu <0$ and $P(\mu) < 0$. Note that $ \cA$ has definite spectrum if and only if there exists $ \mu \in \R$ such that $ J\cA - \mu J >0$. If  there exists $ \mu < 0$  such that $P(\mu) <0$ then $P(\lam)$ is hyperbolic. In any case, if $P(\lam)$ is hyperbolic then $\Lam(\cA) \subset \R$ and each nonzero eigenvalue of $\cA$ is semisimple.

	(b) Let $\mu <0 $ be not an eigenvalue of $P(\lam)$.  Then by (\ref{congru}), we have
    \[
    \imath(J\cA- \mu J) = \imath( - P(\mu)) + (r, 0, 0) =  ( \imath_{-}(P(\mu)), \imath_{+}(P(\mu)), 0) + (r, 0, 0).
    \]
	Hence  $ \imath_{+}(J\cA - \mu J) = \imath_{-}(P(\mu)) + r$ and $ \imath_{-}(J\cA - \mu J) = \imath_{+}(P(\mu))$. Now, if $\mu > \lam^+_{\min}$ then  by  Theorem~\ref{th:defspec}(e), we have
	 \begin{eqnarray*}  \imath_{-}(P(\mu)) + r= \imath_{+}(J\cA- \mu J) &=&  \#\{ \lam \in \Lambda^+(\cA) : \lam > \mu\}+n , \\ \imath_{+}(P(\mu)) = \imath_{-}(J\cA -\mu J) &=& \#\{ \lam \in \Lambda^+(\cA) : \lam < \mu\}.
\end{eqnarray*}
On the other hand, if  $ \mu < \lam^{-}_{\max}$ then by Theorem~\ref{th:defspec}(e), we have that
	 \begin{eqnarray*} \imath_{-}(P(\mu)) + r= \imath_{+}(J\cA- \mu J) &=&  \#\{ \lam \in \Lambda^{-}(\cA) : \lam < \mu\}+r , \\  \imath_{+}(P(\mu))= \imath_{-}(J\cA -\mu J) &=& \#\{ \lam \in \Lambda^{-}(\cA) : \lam > \mu\}.
\end{eqnarray*}
Compare these results with those in  \cite[Proposition~14.3, Corollary~14.4]{Ves11} for the case that $K := GG^*$ is positive definite.

(c) { Assume that $ P(\lam)$ is overdamped and let $ \mu \in \Lambda(P)\setminus\{0\}.$} By Lemma~\ref{lem:lem1}, the map
\[
\mathcal N(P(\mu)) \longrightarrow \mathcal N(S(\mu)),\quad u \longmapsto \left[\begin{matrix}
		\mu u \\ G^*u
	\end{matrix}\right],
\]
is an isomorphism. Hence it follows from (\ref{congru}) that the map
\[ \mathcal N(P(\mu)) \longrightarrow \mathcal N(\cA -\mu I_{n+r}),\;\; u \longmapsto \left[\begin{matrix}
	\mu C^* u \\ G^*u
	\end{matrix}\right],
 \]
is an isomorphism.  Set $ \mathbf{u} := \left[\begin{matrix}
	\mu C^* u \\ G^*u
	\end{matrix}\right]$ and $ \mathbf{v} := \left[\begin{matrix}
	\mu  u \\ G^*u
	\end{matrix}\right]$. Then $ \cA \mathbf{u} = \mu \mathbf{u}$ and $ S(\mu) \mathbf{v} = 0$.
    Set $K:= GG^*.$ Then it follows that
\begin{eqnarray*}
    \mathbf{u}^*J \mathbf{u} &=& -\mu^2 u^*Mu + u^*Ku = - \mathbf{v}^* \frac{\partial }{\partial \mu}S(\mu) \mathbf{v}\\ &=&
	- \mu( 2\mu u^*Mu +  u^*Du) =-\mu u^*\frac{\partial }{\partial \mu}P(\mu)u,
\end{eqnarray*}
which shows that for all $ u \in \mathcal N(P(\mu))$, positive and negative type eigenvalues of the pencil $S(\lam)$ and the polynomial $P(\lam)$ can be defined via   $\mathbf{v}^* \frac{\partial }{\partial \mu}S(\mu) \mathbf{v}$ and $u^*\frac{\partial }{\partial \mu}P(\mu)u$, respectively.
	
 Set $ \Delta(u) := (u^*Du)^2 - 4(u^*Mu) (u^*Ku)$. Since $ \mu^2 u^*Mu + \mu u^*Du+ u^*Ku = 0$, we have { $ \dm{ \mu_{\pm}  = \frac{ - u^*Du \pm \sqrt{\Delta(u)}}{2 u^*Mu}} < 0$.} A simple calculation, see \cite[p.122]{Ves11}, shows  that
\begin{equation}\label{jposevl} \mathbf{u}^*J\mathbf{u} = \frac{ \Delta (u) \pm  u^*Du\sqrt{\Delta(u)} }{2 u^*Mu} = \mp \mu_{\pm} \sqrt{\Delta(u)}.
\end{equation}
As $ \mu_{\pm} <0$, the eigenvector $\mathbf{u}$ of $\cA$ is $J$-positive or $J$-negative according to the sign chosen in the last equality in (\ref{jposevl}). For instance, if we choose
\[
\mu_{+}  = \frac{ - u^*Du + \sqrt{\Delta(u)}}{2 u^*Mu},
\]
then $ \mathbf{u}J\mathbf{u} = - \mu_{+} \sqrt{\Delta(u)} >0$. Hence $ \mu_{+}$ is a $J$-positive eigenvalue of $\cA$.  Similarly, the eigenvalue $\mu_{-}$ is a $J$-negative eigenvalue of $\cA$. Finally, if $ \Delta(u) = 0$ then $\mathbf{u}^*J\mathbf{u} = 0$ which shows that $ \mu = - u^*Du/({2 u^*Mu})$ is a multiple eigenvalue of mixed type. 	
	}
\end{remark}


\section{Deflation and damping of purely imaginary eigenvalues}
\label{Sec:Evolutionofpurimageig}
In this section we discuss the deflation of nonzero purely imaginary eigenvalues and we also discuss the influence of (additional)
damping on these eigenvalues.

Consider the two matrix polynomials $P_0(\lam) := \lam^2 M + K$ and  $P(\lam) := \lam^2 M + \lam D + K$,  where  $M=M^*>0$, $K=K^*\geq 0$, and  $D=D^*\geq 0$. (Note that the following results will also hold if $K$ is indefinite.)
The eigenvalues of $P_0$ are all purely imaginary and typically damping is used to move them from the imaginary axis. In other applications  damping is already active but  this has not yet moved  all eigenvalues  from the imaginary axis. In this case one often uses design changes in the physical system or feedback to move these purely imaginary eigenvalues. { We will show that such additional damping if appropriately chosen can move eigenvalues from the imaginary axis but will not move eigenvalues on the imaginary axis.} Another  application in industrial practice is to compute the eigenvectors of $P_0(\lam)$ and use these for modal truncation in $P(\lam)$ to produce a small size quadratic problem which is assumed to give good approximations to relevant eigenvalues, see \cite{GraMQSW16,MehS11} for an analysis of this approach in the context of brake squeal.  The reduced model is then used for the calculation of extra damping to move the imaginary eigenvalues.

In the following, we introduce a deflation method  for the purely imaginary eigenvalues. 
Consider first the spectral properties of $P(\lambda)$, see \cite{MehMW18} for the matrix pencil case.

\begin{theorem}\label{dfl1}  Consider  a matrix polynomial $ P(\lam) := \lam^2M+\lam D+K$,  where $0<M=M^*\in \mathbb C^{n\times n}$,  $0\leq D=D^*\in \mathbb C^{n\times n}$, and $ 0 \leq K= K^* \in \mathbb{C}^{n\times n}$. Let $ \omega \in \R$ be nonzero.

\begin{itemize}
\item[(a)] Then $i \omega \in \Lam(P)$ if and only if $ \rank\left(\left[\begin{array}{cc} P_0(i\omega) & D\end{array}\right]\right) < n$.  If $i\omega \in \Lam(P)$ then $ i\omega \in \Lam(P_0)$ and $ \mathcal N(P(i\omega)) =\mathcal  N(P_0(i\omega))\cap \mathcal N(D)$.

\item[(b)]Suppose that $i \omega \in \Lam(P)$ and  $\rank\left(\left[\begin{array}{cc} P_0(i\omega) & D\end{array}\right]\right) =n-p$. Then the geometric  multiplicity of $i\omega $ is $p$.     Consider  the QR factorization
\[
Q^*\left[\begin{array}{cc} P_0(i\omega) & D\end{array}\right] = \left[\begin{array}{c} R \\  0_{p, {2n}}\end{array}\right],
\]
where $R$ has full row rank and $Q :=  \left[\begin{array}{cc} Q_1 & Q_2\end{array}\right]$ is unitary with $ Q_2 \in \C^{n\times p}$. Then the columns of $Q_2$ span $\mathcal N(P(i\omega))$. Define $X_2 := Q_2 (Q_2^*MQ_2)^{-1/2}$ and consider the QR factorization
\[
Y^* (MX_2) = \left[\begin{array}{c} R_1 \\  0_{n-p,p}\end{array}\right],
\]
where
    $Y:= \left[\begin{array}{cc} Y_1 & Y_2 \end{array}\right]$ is unitary and $ Y_2  \in \C^{n\times (n-p)}$. Set $ X_1 := Y_2 (Y_2MY_2)^{-1/2}$  and $X := \left[\begin{array}{cc} X_1 & X_2 \end{array}\right]$. Then { $X$ is M-unitary and }
 \[
        X^*P(\lam)X = \diag(
        \widetilde{P}(\lam),  (\lam^2+\omega^2)I_p)
\]
 such that $\Lam(\widetilde{P}) = \Lam(P)\setminus\{\pm i\omega\}$, where $\widetilde{P}(\lam) := X_1^*P(\lam)X_1$.
\end{itemize}
\end{theorem}

\proof First, note that $P(\lam) = P_0(\lam) + \lam D$ and suppose that $ i \omega \in \Lam(P)$. Then there exists a nonzero vector $x$ such that
$P(i\omega) x = (P_0(i\omega) + i\omega D  )x=0$.  Then      $ x^*P_0(i\omega) x+ i \omega x^*Dx =   x^*P(i\omega)x = 0$, which implies that $x^*Dx =0$ as $x^*P_0(i\omega) x$ real and $ \omega \neq 0$. Since  $D$ is semidefinite, we have  $ Dx =0$. Consequently, we have $P_0(i\omega) x = P_0(i \omega) x+ i\omega D x = P(i\omega )x =0$. This shows that $ i \omega \in \Lam(P_0),$ $\mathcal N(P(i\omega))\subset \mathcal N(P_0(i\omega))$ and  $$ \rank\left(\left[\begin{array}{cc} P_0(i\omega) & D\end{array}\right]\right) =  \rank\left(\left[\begin{array}{c} P_0(i\omega) \\  D\end{array}\right]\right)< n.$$

Conversely, suppose that $\rank\left(\left[\begin{array}{cc} P_0(i\omega) & D\end{array}\right]\right) < n$. Then there exists a nonzero vector $u$ such that $ u^*\left[\begin{array}{cc} P_0(i\omega) & D\end{array}\right] = 0$ which implies that $ u^* P_0(i\omega) =0 $ and $ u^*D =0$. Therefore, we have  $P(i\omega)u = P_0(i\omega) u+ i\omega Du =0$ which implies that $ i \omega \in \Lam(P)$ as well as $ i\omega \in \Lam(P_0)$.

 By part (a), we have $ x \in \mathcal N(P(i\omega))$ if and only if $  x \in  \mathcal N(P_0(i\omega)) \cap \mathcal N(D)$ if and only if $ x^* \left[\begin{array}{cc} P_0(i\omega) & D\end{array}\right]=0$. This shows that $p$ is the geometric multiplicity of $i\omega$ as an eigenvalue of $P(\lam)$. Observe that $Q_2^*\left[\begin{array}{cc} P_0(i\omega) & D\end{array}\right] = 0_{p, 2n}$. Hence $P_0(i\omega) Q_2 =0$ and $DQ_2 =0$ which yields $P(i\omega) Q_2 =0$. This shows that the columns of $Q_2$ form an orthonormal basis of $\mathcal N(P(i\omega))$.

By  construction, we have $X^*MX= I_n$.   Note that $ DX_2 =0$  and $ \omega^2 MX_2 = KX_2$ which together with $X^*MX = I_n$ then yields that $ X_1^*P(\lam)X_2 =0$ and $ X_2^*P(\lam) X_2 = (\lam^2+ \omega^2)I_k$.
Since the geometric multiplicity of $i\omega $ is $p$, it follows that $ \pm i\omega \notin \Lam(\widetilde{P})$.
 \eproof

\begin{remark} \label{rem:repeated}{\rm
 Let $ \sig_0 := \{ \pm i \omega_1, \ldots, \pm i \omega_{\ell}\}\subset \Lam(P)$. Suppose that for $j = 1, 2, \ldots, \ell$ the geometric multiplicity of $ i\omega_j$ is $m_j$.
 Set $ m := m_1+\cdots+m_{\ell}$.
 Then repeated application of Theorem~\ref{dfl1}(b) shows that  there exists a nonsingular matrix $X :=\left[\begin{array}{cc} X_1 & X_2\end{array}\right]$ with $ X_2 \in \C^{n\times m}$ such that $  X^*P(\lam) X =  \diag( \widetilde{P}(\lam), \widetilde{P}_0 (\lam))$,
where  $\widetilde{P}_0(\lam) :=  \diag( (\lam^2+\omega_1^2)I_{m_1}, \, \ldots, \, (\lam^2+\omega_{\ell}^2)I_{m_{\ell}})$ and $\widetilde{P}(\lam) := X_1^*P(\lam) X_1$.  Furthermore,  $ \Lam(\widetilde{P}) = \Lam(P)\setminus \sig_0$ and $ \Lam(\widetilde{P}_0) = \sig_0$. This fact is also proved in \cite[Proposition~15.6]{Ves11}.
}
\end{remark}

Theorem~\ref{dfl1} shows the effect of the semidefinite damping matrix on the  purely imaginary eigenvalues of $P_0(\lam)$.  Observe that for any semidefinite damping $D$, we have  $ \Lam(P) \cap i\R \subset \Lam(P_0)\cap i\R$.  This means that the purely imaginary eigenvalues of $P(\lam),$ if any, are the undamped eigenvalues (i.e., the eigenvalues of $P_0(\lam)$)  that remain stationary in $i\R$ after the damping $D$ is applied.

Furthermore, when the damping $D$ is applied, an eigenvalue $i\omega$ (counting multiplicity)  will leave the imaginary axis  if and only if $ \rank\left(\left[\begin{array}{cc} P_0(i\omega) & D\end{array}\right]\right) =n$. Thus, if
$ \Lam(P_0) = \{ \pm i \omega_1, \ldots, \pm i \omega_{\ell}\}$, then
all the undamped eigenvalues  will leave the imaginary axis when the damping $D$ is applied if and only if \\ $ \rank\left(\left[\begin{array}{cc} P_0(i\omega_j) & D\end{array}\right]\right) =n$ for $ j= 1, 2, \ldots, \ell$. In particular, if the damping matrix $D$ is definite then all undamped eigenvalues will leave the imaginary axis. This also shows the well-known fact that the eigenvalues of $P(\lam)$ have negative real part when $M$, $D$, and $K$ are all positive definite.

We now analyze in more detail how to construct damping matrices $D$ that remove undamped eigenvalues from the imaginary axis.   Consider $P_0(\lam) := \lam^2 M + K$ with $M>0$ and $K>0$. Then for any positive semidefinite damping matrix $D$ and $ P(\lam) := P_0(\lam) +\lam D$, we have  $ \Lam(P)\cap i\R \subset  \Lam(P_0)$.
More precisely, let
  $ i \omega \in \Lam(P_0)$ with $\omega  \neq 0$, then by Theorem~\ref{dfl1} there are three possible cases:
(a)
$\pm i\omega \notin \Lam(P)$  or
(b) $\pm i\omega \in \Lam(P)$ and $ \mathcal N(P(i\omega)) \subsetneqq \mathcal N(P_0(i\omega))$
or
 (c)
 $\pm i\omega \in \Lam(P)$ and $\mathcal  N(P(i\omega)) = \mathcal N(P_0(i\omega))$.

In case (a), the damping $D$ removes $\pm i\omega$ from the imaginary axis. In case (b), $i\omega$ is a multiple eigenvalue and the damping $D$ removes a few copies of the eigenvalue $\pm i\omega$ (counting multiplicity) from the imaginary axis and a few copies of  $\pm i\omega$ remain purely imaginary with reduced geometric multiplicity. Finally, in case (c), $\pm i\omega$  is completely unaffected by the damping $D$.

 Our next result shows that if $i\omega$ is an undamped eigenvalue of geometric multiplicity $p$ then a semidefinite damping matrix that  removes  $ i\omega$ completely from the imaginary axis must have rank  at least $p$.

\begin{theorem} \label{dfl3}
 Consider  a matrix polynomial $ P(\lam) := \lam^2M+\lam D+K$, where $0<M=M^*\in \mathbb C^{n\times n}$,  $0\leq D=D^*\in \mathbb C^{n\times n}$, and $ 0 \leq K = K^* \in \mathbb C^{n\times n}$. Let $ P_0(\lam) := \lam^2 M + K$ and let $ i\omega \in \Lam(P)$ with  geometric multiplicity $p$. Consider the QR factorization $ Q^*\left[\begin{array}{cc} P_0(i\omega) & D\end{array}\right] =  \left[\begin{array}{c} R \\  0_{p, {2n}}\end{array}\right]$, where $R$ has full row rank and $Q :=  \left[\begin{array}{cc} Q_1 & Q_2\end{array}\right]$ is unitary with $ Q_2 \in \C^{n\times   p}.$  Then the columns of $Q_2$ span $\mathcal N(P(i\omega))$.

 Let $ 0 < T = T^* \in \C^{p,p}$. Consider the family of matrix polynomials  $ \widehat{P}_{t}(\lam) := P(\lam)+\lam \, t MQ_2TQ^*_2M$ for all $t>0$.  Then there exists { a $M$-unitary}  matrix $X = \left[\begin{array}{cc} X_1 & X_2 \end{array}\right]$ with  $X_2 := Q_2 (Q_2^*MQ_2)^{-1/2}$  such that for all $ t >0$,
\begin{eqnarray*}
    X^*P(\lam)X &=& \diag(
    \widetilde{P}(\lam), (\lam^2+\omega^2)I_p),\\
X^*\widehat{P}_t(\lam)X &=& \diag( \widetilde{P}(\lam), \lam^2 I_p+ \lam\, t\widehat{T} +\omega^2I_p),
\end{eqnarray*}
 where $\widehat{T} := (Q_2^*MQ_2)^{1/2}T(Q_2^*MQ_2)^{1/2}$ and $\widetilde{P}(\lam) := X_1^*P(\lam)X_1$.
 Furthermore, we have $\Lam(\widetilde{P}) = \Lam(P)\setminus\{\pm i\omega\}$ and the eigenvalues of $\lam^2I_p+ \lam\, t\widehat{T} +\omega^2I_p$ have negative real parts for all $t>0$. Thus, in the  damped system $ \widehat{P}_{t}(\lam)$, for all $ t >0$, $\pm i\omega \notin \Lam(\widehat{P}_{t})$.

Let  $W$ be an $n\times m$ matrix such that $\rank(W) = m <p$. Then for any $t>0$ the damped  system { $ \widetilde{P}_{t}(\lam) := P(\lam)+\lam \, t WW^*$
 has  undamped eigenvalues  $\pm i\omega $, i.e., $\pm i\omega \in \Lam(\widetilde{P}_{t})$ for all $ t >0$.}
  \end{theorem}

\proof Set $ \widehat{D} := MQ_2 TQ_2^*M = MX_2 \widehat{T}X_2^*M $.  Note that the columns of  $X_2$ form an $M$-orthonormal  basis of $\mathcal N(P(i\omega))$. For the given $X_2$ construct $X_1$ as in Theorem~\ref{dfl1}(b), so that by construction $X^*MX= I_n$.  We have  $ DX_2 =0$  and $ \omega^2 MX_2 = KX_2$ which together with $X^*MX = I_n$ yields $ X_1^*P(\lam)X_2 =0$ and $ X_2^*P(\lam) X_2 = (\lam^2+ \omega^2)I_p$. Next, observe that $\widehat{D} X_1 =0 $ and $ X_2^*\widehat{D}X_2 = \widehat{T}$. Hence by Theorem~\ref{dfl1}, for all $ t>0$ we have
\begin{eqnarray*}X^*P(\lam)X &=& \diag( \widetilde{P}(\lam), (\lam^2+\omega^2)I_p), \\
X^*\widehat{P}_t(\lam)X& =& \diag(
\widetilde{P}(\lam), \lam^2I_p+ \lam\, t\widehat{T} +\omega^2I_p)
\end{eqnarray*}
and $\Lam(\widetilde{P}) = \Lam(P)\setminus\{\pm i\omega\}$. Since $\widehat{T}$ is positive definite, again by Theorem~\ref{dfl1}, the eigenvalues of  $\lam^2I_p+ \lam\, t\widehat{T} +\omega^2I_p$ lie in the open left half complex plane.

Finally, set $\widehat{D}_t := D +t WW^*$. Then $ W = Q_1Y_1+Q_2Y_2$ for some $(n-p)\times m$ matrix $Y_1$ and $p\times m$ matrix $Y_2$. Since $ m<p$, there exists a nonzero vector $u$ such that $ u^*Y_2 =0$. With $ v := Q_2u$, then $v^*W = u^*Y_2 = 0$. Consequently,  $v^*\widehat{D}_t = v^*(D +t WW^*) = 0$ and $ P_0(i\omega) v=0$ which shows that
$ \rank\left(\left[\begin{array}{cc} P_0(i\omega) & \widehat{D}_t\end{array}\right]\right) < n$ for all $t >0$.  Hence by Theorem~\ref{dfl1}, we have $ \pm i\omega \in \Lam(\widehat{P}_{t})$ for all $t>0$.
\eproof


\begin{remark} \label{rem1}{\rm
Observe that for all $t>0$ the damped system $\widehat{P}_t(\lam)$ in Theorem~\ref{dfl3} removes $\pm i\omega $ from the imaginary axis and leaves the remaining spectrum of $P(\lam)$ (including the partial multiplicities of eigenvalues)  completely unaffected. The $2p$ eigenvalues $\pm i\omega$ (counting multiplicity) evolve as eigenvalues of
$\Theta_t(\lam) := \lam^2I_p+ \lam\, t\widehat{T} +\omega^2I_p$ for $ t >0$. The complex eigenvalues of $ \Theta_t(\lam)$, if there are any, lie on the semicircle $|\lam|^2 = \omega^2$ in the left half complex plane. Indeed, if $u$ is an eigenvector with $u^*u=1$ then $\lam^2 + \lam t u^* \widehat{T} u + \omega^2 = 0$ implies that
\[\lam_{\pm}(t) = \frac{ - t u^*\widehat{T}u \pm \sqrt{  (t u^*\widehat{T}u)^2 - 4 \omega^2}}{2}.
\]
For complex eigenvalues, we have $(t u^*\widehat{T}u)^2 - 4 \omega^2 <0$ which shows that $ |\lam_{\pm}(t)|^2 = \omega^2$.

If  $(t u^*\widehat{T}u)^2 - 4 \omega^2 >0$ then we have two distinct real eigenvalues given by
\begin{eqnarray*}
\lam_{+}(t) &=& \frac{ - 2 \omega^2}{  t u^*\widehat{T}u +\sqrt{  (t u^*\widehat{T}u)^2 - 4 \omega^2}},\ \lam_{-}(t) =  -\frac{  t u^*\widehat{T}u + \sqrt{  (t u^*\widehat{T}u)^2 - 4 \omega^2}}{2}.
\end{eqnarray*}
It follows from (\ref{jposevl}) and the discussion there that $ \lam_{+}(t)$ is an eigenvalue of positive type and $\lam_{-}(t)$ is  an eigenvalue of negative type. Also $ \lam_{+}(t) \rightarrow 0$ and  $ \lam_{-}(t) \rightarrow -\infty$ as $ t \rightarrow \infty$.

If $t  \lam_{\min}(\widehat{T})  > 2 \omega$ then $  (t\, x^*\widehat{T}x)^2 - 4 \omega^2 (x^*x)^2 >0 $ for all $x \neq 0$. Hence, for $ t > 2\omega /{\lam_{\min}(\widehat{T})} $, $\Theta_t(\lam)$ is hyperbolic and has $2p$ real eigenvalues (counting multiplicities) consisting of $p$ eigenvalues of positive type and $p$ eigenvalues of negative type. Furthermore,  the $p$ positive type eigenvalues of $\Theta_t(\lam)$ converge to $0$  and the $p$ negative type eigenvalues converge to $-\infty$ as   $t \rar \infty$. It also  follows that the eigenvalues $\pm i\omega$, once they leave the imaginary axis, can never return to the imaginary axis as nonzero/noninfinite purely imaginary eigenvalues of $\widehat{P}_t$ for any $ t>0$.

Finally, if $(t u^*\widehat{T}u)^2 - 4 \omega^2 = 0$ then $ \lam_{\pm}(t) = - t u^* \widehat{T} u/2$ is a multiple eigenvalue of mixed type.
}
\end{remark}


\begin{remark} \label{rem3} {\rm
{ Let $ X =\left[\begin{array}{cc} X_1 & X_2 \end{array}\right]$ with $X_2 := Q_2 (Q_2^*MQ_2)^{-1/2}$  be as in Theorem~\ref{dfl3}.} Now, replace the damping matrix $ MQ_2TQ^*_2M$  in Theorem~\ref{dfl3} with $ MX_2TX^*_2M$, i.e., for  $t>0$, consider the damped system $ \widehat{P}_{t}(\lam) := P(\lam)+\lam \, t MX_2TX^*_2M$. {   Since $X$ is $M$-unitary, }for all $ t >0$ we have
\begin{eqnarray*}X^*P(\lam)X &=& \diag( \widetilde{P}(\lam), (\lam^2+\omega^2)I_p), \\
X^*\widehat{P}_t(\lam)X& =& \diag( \widetilde{P}(\lam),\lam^2I_p+ \lam\, tT +\omega^2I_p),
\end{eqnarray*}
 where  $\widetilde{P}(\lam) := X_1^*P(\lam)X_1$ and  $\Lam(\widetilde{P}) = \Lam(P)\setminus\{\pm i\omega\}$.  The  eigenvalues of $\Theta_t(\lam) :=\lam^2I_p+ \lam\, t T +\omega^2I_p$ have nonzero negative real parts for all $t>0$. Thus, the damping $t MX_2TX_2^*M$  removes $\pm i\omega $ completely from the imaginary axis and the eigenvalues $\pm i\omega$ move to the eigenvalues of $\Theta_t(\lam)$ for all $ t>0$.  In particular, choosing $ T := \diag(d_1, \ldots, d_p) > 0$, we have $ \Theta_t(\lam) =  \diag(q_1(\lam, t), \ldots, q_p(\lam, t))$, where $ q_j(\lam, t) := \lam^2 + \lam t d_j + \omega^2$ for $ j = 1,\ldots,p$. Therefore, by choosing $d_j$ appropriately, the eigenvalues $\pm i\omega$ can be moved to  any pre-specified locations in the complex plane without changing the remaining eigenvalues.
 }
 \end{remark}

It turns out that if the eigenvalues of $P_0(\lam)$ are simple then a  semidefinite damping matrix of rank one can remove all the  eigenvalues of $P_0(\lam)$ from the imaginary axis.

\begin{proposition} Let   $ P_0(\lam) := \lam^2 M+ K$,  where $ 0 < M = M^*,\  0 <K = K^* \in \C^{n\times n}$. Let $X$ be $M$-unitary  such that $X^*KX = \diag(\omega_1^2, \ldots, \omega_n^2)$.  Set $ v := MXu$, where $ u :=[u_1, \ldots, u_n]^\top \in \C^n $ and  $ u_j \neq 0$ for $ j = 1,\ldots,n$. Consider $ P_t(\lam) := P_0(\lam) + \lam \, t vv^*$. If all the eigenvalues of $P_0(\lam)$ are simple then, for all $t>0$, the eigenvalues of $P_t(\lam)$  have negative real part.
\end{proposition}
\proof Note that the columns of $X$ are eigenvectors of $ P_0(\lam)$. Let $x_j $ be the $j$-th column of $X$, that is, $ x_j := Xe_j$. Then $ P_0(i\omega_j)x_j =0$ for $ j = 1,\ldots,n$.
Since $ i\omega_j$ is simple, $ \pm i \omega_j \in \Lam(P_t)$ if and only if $ t vv^* x_j = 0$ if and only if $v^*x_j=0$. Now $ v^*x_j = u^* X^*Mx_j = u^*X^*MXe_j = u_j \neq 0$. This shows that $ t vv^*x_j \neq 0$ for $j = 1,\ldots,n$ and $ t>0$. Since $ \Lam(P_t) \cap i\R \subset \Lam(P_0)$ for all $ t >0$, we conclude that  all eigenvalues of $ P_t(\lam)$ have negative real part.
\eproof

 The damping matrix $\widehat{D} := M Q_2 TQ_2^*M $ in Theorem~\ref{dfl3} removes the imaginary eigenvalues $\pm i\omega$ of $P(\lam)$ from the imaginary axis leaving the other eigenvalues of $P(\lam)$  unchanged. On the other hand, the damping matrix $\widetilde{D} := Q_2TQ_2^*$ can also be used to remove the eigenvalues $\pm i\omega$ from the imaginary axis even though a decomposition of $\widehat{P}_t(\lam) := P(\lam) + \lam \, t \widehat{D}$ as given in Theorem~\ref{dfl3} may not be possible. 

\begin{proposition}  \label{new:prop2} Consider $P(\lam) := \lam^2 M + \lam D + K$ and  $P_0(\lam) = \lam^2M+ K$,  where $0 < M = M^*,\ 0 \leq D = D^*,\  0 \leq K = K^* \in \C^{n\times n}$. Let $ i\omega \in \Lam(P)$ have  geometric multiplicity $p$ and let $ Q_2$ and $T$ be as in Theorem~\ref{dfl3}.  Consider  $ \widehat{P}_{t}(\lam) := P(\lam)+\lam \, t Q_2TQ^*_2$.  Then  for all $ t>0$, $ \pm i\omega \notin \Lam(\widehat{P}_{t})$.
  \end{proposition}

\proof
Set $\widehat{D}_t := D+t Q_2TQ_2^*.$
It is easy to see that $\mathcal N(\widehat{P}_t(i\omega) ) \subset \mathcal N(P(i\omega))$ for $ t >0$. Indeed, by Theorem~\ref{dfl1}, we have $\mathcal N(\widehat{P}_t(i\omega) ) = \mathcal N(P_0(i\omega))\cap \mathcal N(\widehat{D}_t)$. Since $D\geq 0$ and $T >0$, we have  that $ x \in \mathcal N(\widehat{D}_t)$ if and only if $Dx + t(Q_2TQ_2^*)x = 0$ if and only if $ Dx = 0$  and $Q_2TQ_2^* x =0$. Hence $\mathcal N(\widehat{P}_t(i\omega) ) \subset \mathcal N(P(i\omega))$ for $ t >0.$

Note that the columns of $Q_2$ form an orthonormal basis of $\mathcal N(P(i\omega))$ and $DQ_2 =0$. Let $ x \in \mathcal  N(P_0(i\omega))\cap \mathcal N(\widehat{D}_t)$. Then $x = Q_2y$ for some $ y \in \C^p$. Now  $\widehat{D}_tx = DQ_2y + t Q_2Ty = tQ_2Ty$. Hence $\widehat{D}_tx =0 $ implies that $Q_2Ty = 0$ and thus $ y = 0$ and $ x = Q_2 y = 0$. This shows that $ x^*\left[\begin{array}{cc} P_0(i\omega) & \widehat{D}_t\end{array}\right]  = 0 $ implies that $ x = 0$. Hence, we have that $ \rank\left(\left[\begin{array}{cc} P_0(i\omega) & \widehat{D}_t\end{array}\right]\right) = n$ for  all $t>0$. Thus, by Theorem~\ref{dfl1}, we have  $ \pm i \omega \notin \Lam(\widehat{P}_{t})$ for all $t >0.$
\eproof

\begin{remark}\label{rem2}{\rm
 Suppose that $ i \omega \in \Lam(P)$.  Consider $ \widehat{P} (\lam) := P(\lam) + \lam \widehat{D}$, where $\widehat{D }\geq 0$ is such that
 $\rank\left(\left[\begin{array}{cc} P_0(i\omega) & \widehat{D}\end{array}\right]\right) =n$. By Theorem~\ref{dfl1}, we have  $ \pm i\omega \notin \Lam(\widehat{P})$. If $F=F^* \geq 0$ then it can be shown that $\pm i\omega$ are not the eigenvalues of $\widehat{ P}(\lam) + \lam F = P(\lam) + \lam (\widehat{D}+F)$. Indeed, this follows from the fact that
  $\rank\left(\left[\begin{array}{cc} P_0(i\omega) & \widehat{D} +F\end{array}\right]\right) =n$ for any $ F \geq 0$. This means that if purely imaginary eigenvalues leave the imaginary axis under a positive semidefinite damping then these eigenvalues remain away from the imaginary axis for any subsequent further positive semidefinite damping.    Therefore, Theorem~\ref{dfl3} 
  can be used to remove one pair of imaginary eigenvalues $\pm i\omega$ of $P(\lam)$ at a time.
  }
\end{remark}

 We now summarize the evolution of purely imaginary eigenvalues of $P(\lam)$ under parameter dependent positive semidefinite damping $tD$ for $t >0$. \\

\begin{corollary} \label{eig:imag1} Let   $ P_0(\lam) := \lam^2 M+ K$ and $ P_t(\lam) := \lam^2M+\lam tD+ K$ for $ t>0$, where $ M>0$, $D\geq 0$, and $ K >0$. Then $\Lam(P_t)\cap i\R \subset \Lam(P_0)$ for all $t >0$.

(a) { Let $ i \omega \in \Lam(P_0)$ and let $m$ be the geometric multiplicity of $i\omega.$ } If $\rank\left(\left[\begin{array}{cc} P_0(i\omega) & D\end{array}\right]\right) = n$ then $ \pm i \omega \notin \Lam(P_t)$ for all $t>0$.
On the other hand, if $\rank\left(\left[\begin{array}{cc} P_0(i\omega) & D\end{array}\right]\right) = n-p$ then $\pm i \omega \in \Lam(P_t)$ for all $t>0$. Further, $p$ is the geometric multiplicity of $\pm i\omega$ as an eigenvalue of $P_t(\lam)$   and $ \mathcal N(P_t(i\omega)) = \mathcal N(P_1(i\omega))$ for all $t>0$. Hence $P_t(\lam)$ removes $(m-p)$ eigenvalues $\pm i\omega$ (counting multiplicity) from the imaginary axis.

(b) Let $ \Lam(P_1)\cap i\R = \{ \pm i \omega_1, \ldots, \pm i \omega_{\ell}\}=:\sig_0$ and $\rank\left(\left[\begin{array}{cc} P_0(i\omega_j) & D\end{array}\right]\right) = n-p_j$ for $j = 1, 2, \ldots, \ell$.  Set $ p := p_1+\cdots+p_{\ell}$. Then there exists { a $M$-unitary} matrix $X :=\left[\begin{array}{cc} X_1 & X_2\end{array}\right]$ with $ X_2 \in \C^{n\times p}$ such that $  X^*P_t(\lam) X =  diag( \widetilde{P}_t(\lam), \widetilde{P}_0 (\lam))  \, \mbox{ for all }\, t >0,$  where  $ \widetilde{P}_0(\lam) :=  \diag( (\lam^2+\omega_1^2)I_{p_1}, \, \ldots, \, (\lam^2+\omega_{\ell}^2)I_{p_{\ell}})$ and $\widetilde{P}_t(\lam) := X_1^*P_t(\lam) X_1$.
Furthermore, all eigenvalues of $\widetilde{P}_t(\lam)$ have negative real parts, $ \Lam(\widetilde{P}_t) = \Lam(P_t)\setminus \sig_0$ and $ \Lam(\widetilde{P}_0) = \sig_0$ for all $ t >0$.  The geometric multiplicity of a purely imaginary eigenvalue of $P_t(\lam)$ remains the same  for all $ t>0$.

\end{corollary}

\proof Note that $\rank\left(\left[\begin{array}{cc} P_0(i\omega) & tD\end{array}\right]\right) = \rank\left(\left[\begin{array}{cc} P_0(i\omega) & D\end{array}\right]\right)$ for all $ t>0$. Hence the conclusions in (a) follow. The matrix $X$ as constructed in Theorem~\ref{dfl1} depends on an $M$-orthonormal basis of $ \mathcal N(P_1(i\omega_j))$. Since $ \mathcal N(P_1(i\omega_j))  =  \mathcal N(P_t(i\omega_j))$ for all $t >0$, the assertions in (b) follow from Theorem~\ref{dfl1}.
\eproof

\subsection{Numerical methods for the deflation of purely imaginary eigenvalues.} \label{num}
{ We now present two algorithms that implement Theorem~\ref{dfl3} to deflate a purely imaginary eigenvalue. These are based on rank-revealing factorizations and we specifically discuss the rank-revealing QR and singular value decomposition.} Let $  0< M,\  0 \leq K,\  D= D \in \C^{n\times n}$ with $D$ semidefinite.  Consider  $P_0(\lam) := \lam^2 M + K $ and $ P(\lam) := \lam^2 M + \lam D + K.$  We describe two numerical methods for deflating the purely imaginary eigenvalues of $ P(\lam)$. The first method is based on an $M$-unitary QR factorization, {which is actually a version of the Cholesky QR decomposition, see e.g. \cite{StaW02}.}
		
Let $ i\omega  \in \Lambda(P) $ and $\omega \neq 0$. For  $m : = \rank( \left[\begin{array}{cc} P_0(i\omega) & D \end{array}\right]) < n$, set $ p := n-m$. Consider a rank revealing  QR factorization
\[
\left[\begin{array}{cc} P_0(i\omega) & D \end{array}\right] \Pi =  Q  \left[\begin{array}{cc} R_{11} & R_{12} \\ 0 & 0 \end{array}\right],
\]
where $   Q^*Q = I_n$,  $ R_{11} \in \C^{m \times m}$ is nonsingular and $\Pi \in \C^{2n \times  2n}$ is a permutation matrix. Partition $ Q=\left[\begin{array}{cc} Q_1 & Q_2 \end{array}\right]$, where $ Q_1 \in \C^{n\times  m}$ are the leading columns of $Q$. Then it follows that $ Q_2^* \left[\begin{array}{cc} P_0(i\omega) & D \end{array}\right] = 0$ which shows that   $\mathrm{span}(Q_2) = \mathcal{N}(D) \cap \mathcal{N}(P_0(i\omega))$.
	
Next, compute a Cholesky factorization $Q_2^* M Q_2=C_2 C_2^*$ and define
$X_2:= Q_2(C_2^*)^{-1}$ by solving the linear system $X_2 C_2^*=Q_2$ for $X_2$.
Then we have $X_2^* M X_2= C_2^{-1} Q_2^* M Q_2 C_2^{-*} = C_2^{-1}  C_2 C_2^*C_2^{-*}=I_p$ and $\mathrm{span}(X_2) = \mathrm{span}(Q_2)  = \mathcal{N}(D) \cap \mathcal{N}(P_0(i\omega))$.   Note that  $  P_0(i\omega) X_2 = 0$ which implies that $ \omega^2 M X_2 = KX_2$ which in turn shows that $ X_2^* KX_2 = \omega^2 I_p$. Since $DX_2 = 0$, we have
\[
X_2^*P(\lam) X_2 = \lambda ^2 X_2^*MX_2 + \lam X_2^*DX_2 + X_2^*KX_2 = (\lam^2 + \omega^2)I_p.
\]
Next, compute a QR factorization $ MX_2 = U\left[\begin{array}{c} R \\ 0 \end{array}\right]$, where   $ R \in \C^{p \times p}$ is nonsingular.   Partition $ U =\left[\begin{array}{cc} U_1 & U_2 \end{array}\right]$, where $ U_1 \in \C^{n\times  p}$.
Compute a Cholesky factorization $ U_2^*MU_2 = C_3C_3^*$ and define $ X_1 := U_2 (C_3^*)^{-1}$ by solving the linear system $ X_1C_3^* = U_2$ for $X_1$.  Then $ X_1^* MX_1 = C_3^{-1} U_2^*MU_2 (C_3^*)^{-1} = C_3^{-1}C_3C_3^* (C_3^*)^{-1}=I_m$ and $ X_1^*MX_2 = 0$ as { $ U_2^* MX_2 = 0.$}	Since  $ \omega^2 MX_2 = KX_2$, it follows that $ X_1^* KX_2 =0$.

Then, setting $ X :=  \left[\begin{array}{cc} X_1 & X_2 \end{array}\right]$, by construction we have $ X^* MX = I_n$ and
\[
X^* P(\lam) X =  \left[\begin{array}{cc}  X_1^*P(\lam) X_1  &  X_1^*P(\lam) X_2 \\ X_2^*P(\lam) X_1 & X_2^*P(\lam)X_2 \end{array}\right] = \left[\begin{array}{cc}  X_1^*P(\lam) X_1  &  0 \\ 0 & (\lam^2 +\omega^2) I_p \end{array}\right].
\]
Since $p$ is the geometric multiplicity of $i\omega,$  it follows that $ \pm i\omega$ is not an eigenvalue  of $  {\widetilde{P}(\lam)} := X_1^* P(\lam)X_1$.
	
	This  leads to the following algorithm; for a {\sc Matlab} code, see \cite{AlaMT25_arxiv}.

	
\begin{algorithm} \label{alg3} QR-based method for deflating nonzero purely imaginary eigenvalues of $P(\lambda)$.
	\begin{itemize}
	\item[1.] Compute $m : = \rank( \left[\begin{array}{cc} P_0(i\omega) & D \end{array}\right])$  and set $ p := n-m$.  IF $ p=0$ then STOP as $\pm i\omega$ are not eigenvalues of $P(\lam)$  ELSE proceed as follows.
	\item[2.]  Compute a rank revealing QR factorization
    \[ \left[\begin{array}{cc} P_0(i\omega) & D \end{array}\right] \Pi = Q  \left[\begin{array}{cc} R_{11} & R_{12} \\ 0 & 0 \end{array}\right],
    \]
    where $ R_{11} \in \C^{m \times m}$ is nonsingular and $\Pi \in \C^{2n \times 2n}$ is a permutation matrix.
    \item[3.] Partition $ Q =\left[\begin{array}{cc} Q_1 & Q_2 \end{array}\right]$, where $ Q_1 \in \C^{n\times  m}$ and compute the Cholesky factorization $ Q_2^* M Q_2 = C_2C_2^*$.
	\item[4.] Solve the linear system $X_2 C_2^*=Q_2$ for $X_2$.
	\item[5.] Compute a QR factorization $ MX_2 = U\left[\begin{array}{c} R \\ 0 \end{array}\right]$, where $ R \in \C^{p \times p}$ is nonsingular.
	\item[6.]  Partition $ U =\left[\begin{array}{cc} U_1 & U_2 \end{array}\right]$, where $ U_1 \in \C^{n\times  p}$ and compute the Cholesky factorization $ U_2^*MU_2 = C_3C_3^*$.
    \item [7.] Solve the system $ X_1C_3^* = U_2$ and define $X := \left[\begin{array}{cc} X_1 & X_2 \end{array}\right]$.
	\item[8.] 	Then $ X^* MX = I_n$ and  $ X^* P(\lam) X =  \diag(  X_1^*P(\lam) X_1, (\lam^2 +\omega^2) I_p )$.
	\end{itemize}
\end{algorithm}	
	
The second method is based on SVD.
	%
    %
%
Let $ i\omega  \in \Lambda(P) $  and $\omega\neq 0$. Suppose that   $m : = \rank( \left[\begin{array}{cc} P_0(i\omega) & D \end{array}\right]) < n$ and set $ p := n-m$. Consider the SVD
\[
\left[\begin{array}{cc} P_0(i\omega) & D \end{array}\right] = U \left[\begin{array}{cc} \diag(\sigma_1, \ldots, \sigma_m)  &  0 \\ 0 & 0 \end{array}\right]V^* = : U\Sigma V^*.
\]
With the partition $ U =\left[\begin{array}{cc} U_1 & U_2 \end{array}\right]$, where $ U_1 \in \C^{n\times  m}$, then it follows that $ U_2^* \left[\begin{array}{cc} P_0(i\omega) & D \end{array}\right] = 0$ which shows that   $\mathrm{span}(U_2) = \mathcal{N}(D) \cap \mathcal{N}(P_0(i\omega))$.

Next, compute the Cholesky factorization $ U_2^*MU_2 = C_2 C_2^*$ and define $ X_2 := U_2 (C_2^*)^{-1}$ by solving the linear system $ X_2 C_2^* = U_2$ for $X_2$.   Then we have  that $ X_2^* MX_2 = C_2^{-1} U_2^*MU_2 (C_2^*)^{-1} = C_2^{-1}C_2C_2^* (C_2^*)^{-1}= I_p$  and $\mathrm{span}(X_2) = \mathrm{span}(U_2)  = \mathcal{N}(D) \cap \mathcal{N}(P_0(i\omega))$.   Note that  $  P_0(i\omega) X_2 = 0$ which implies that $ \omega^2 M X_2 = KX_2$ which in turn shows that $ X_2^* KX_2 = \omega^2 I_p$. Since $DX_2 = 0$, we have
\[
X_2^*P(\lam) X_2 = \lambda ^2 X_2^*MX_2 + \lam X_2^*DX_2 + X_2^*KX_2 = (\lam^2 + \omega^2)I_p.
\]
Next, consider  the SVD $ MX_2 = Y\left[\begin{array}{c} \diag(\tau_1, \ldots, \tau_p) \\ 0 \end{array}\right] W^*$.    Partition $ Y =\left[\begin{array}{cc} Y_1 & Y_2 \end{array}\right]$, where $ Y_1 \in \C^{n\times  p}$. Compute the Cholesky factorization $ Y_2^*MY_2 =C_3C_3^*$ and define $ X_1 := Y_2 (C_3^*)^{-1}$ by solving the linear system $ X_1C_3^* = Y_2$ for $X_1$.  Then $ X_1^* MX_1 = C_3^{-1} Y_2^*MY_2 (C_3^*)^{-1} = C_3^{-1} C_3C_3^* (C_3^*)^{-1}= I_m$ and $ X_1^*MX_2 = 0$.	Since  $ \omega^2 MX_2 = KX_2$, it follows that $ X_1^* KX_2 =0$.

 Define $ X :=  \left[\begin{array}{cc} X_1 & X_2 \end{array}\right]$. Then by construction we have $ X^* MX = I_n$ and
 \[
 X^* P(\lam) X =  \left[\begin{array}{cc}  X_1^*P(\lam) X_1  &  X_1^*P(\lam) X_2 \\ X_2^*P(\lam) X_1 & X_2^*P(\lam)X_2 \end{array}\right] = \left[\begin{array}{cc}  X_1^*P(\lam) X_1  &  0 \\ 0 & (\lam^2 +\omega^2) I_p \end{array}\right].
 \]
As $p$ is the geometric multiplicity of $i\omega,$  it follows that $ \pm i\omega$ is not an eigenvalue  of $ \widehat{P}(\lam) := X_1^* P(\lam)X_1$.

\begin{algorithm} \label{alg4} SVD-based method for deflating purely imaginary eigenvalues of $P(\lambda)$.
\begin{itemize}
		\item[1.] Compute $m : = \rank( \left[\begin{array}{cc} P_0(i\omega) & D \end{array}\right])$ and set $ p := n-m$.  IF $ p=0$ then STOP as $\pm i\omega$ are not eigenvalues of $P(\lam)$  ELSE proceed as follows.
	\item[2.]  Compute the SVD $ \left[\begin{array}{cc} P_0(i\omega) & D \end{array}\right]  = U\Sigma V^*$, where $ U \in \C^{n\times n}$ and $V \in \C^{2n \times 2n}$ are unitary matrices.
	\item[3.] Partition $ U =\left[\begin{array}{cc} U_1 & U_2 \end{array}\right]$, where $ U_1 \in \C^{n\times  m},$  and compute the Cholesky factorization $ U_2^*MU_2 = C_2C_2^*$.
	\item[4.] Solve the linear system $ X_2 C_2^* = U_2$ for $X_2$.
	\item[5.] Compute the SVD  $ MX_2 = Y\left[\begin{array}{c} \diag(\tau_1, \ldots, \tau_p) \\ 0 \end{array}\right] W^*$, where $ Y \in \C^{n\times  n}$ and $ W \in \C^{p, p}$ are unitary.
	\item[6.]  Partition $ Y =\left[\begin{array}{cc} Y_1 & Y_2 \end{array}\right]$, where $ Y_1 \in \C^{n\times  p}$, and compute the Cholesky factorization $ Y_2^*MY_2 = C_3C_3^*.$
	\item[7.] Solve the linear system $ X_1 C_3^* = Y_2$ for $ X_1$ and define $X := \left[\begin{array}{cc} X_1 & X_2 \end{array}\right]$.
		\item[8.] 	Then $ X^* MX = I_n$ and  $ X^* P(\lam) X =  \diag(X_1^*P(\lam) X_1  , (\lam^2 +\omega^2) I_p ).$\end{itemize}
\end{algorithm}

A {\sc Matlab} code for this algorithm is also presented in \cite{AlaMT25_arxiv}.

The matrix  $X := \left[\begin{array}{cc} X_1 & X_2 \end{array}\right]$ computed either by the QR based method or the SVD based method can be used in Remark~\ref{rem3} to construct a damping matrix $\widehat{D} := M X_2 T X_2^* M$  so as to move the eigenvalues $\pm i\omega$ to  pre-specified locations in the complex plane while leaving the other eigenvalues unchanged.

\subsection{Numerical examples}\label{sec:num2}
{ To illustrate our results, in this section we consider two numerical examples.
\begin{example}\label{ex:constructued}{\rm
Consider a problem with
$ n = 10 $,  $    M = I_n$, $K=K_0+K_1$, with \[
K_0 = X_0 \, \mathrm{diag}(25, 49) \, X_0^T,
\]
%
with a real orthonormal matrix $ X_0 \in \mathbb{R}^{n \times p} $, and
$K_1 = G G^T$,
where $ G \in \mathbb{R}^{n \times (r - 2)} $ and the columns of $ G $ are orthogonal to those of $ X_0 $. Then a damping matrix is added so that still has the purely imaginary eigenvalues  $\pm 7i $ and $ \pm 5i $. See the {\sc Matlab} code for this example in \cite{AlaMT25_arxiv}.
computed eigenvalues by   \texttt{polyeig(K, D, M)}  are
\[
\{ 2.9424 \times 10^{-15} + 7i,\ 2.9424 \times 10^{-15} - 7i,\ -2.5853 \times 10^{-15} + 5i,\ -2.5853 \times 10^{-15} - 5i \},
\]
along with two further (numerically zero) eigenvalues $ 1.9334 \times 10^{-15} + 0i$, $-1.1545 \times 10^{-16} + 0i$.

To deflate the pair $  \pm i\omega_1 = \pm 5i $, we apply Algorithm~\ref{alg3}  and obtain a reduced system $ \widehat{M}, \widehat{D}, \widehat{K} $, with eigenvalues
\[
\begin{aligned}
&-3.27\!\times\!10^{-15} \pm 7i,\
-0.63,\ -0.43,
-0.13 \pm 0.99i,\
-0.19 \pm 0.97i,\
-0.49 \pm 0.87i,\\
&-0.44 \pm 0.90i,\
-0.29 \pm 0.95i,\
-0.33 \pm 0.94i,\ 1.49\!\times\!10^{-15},\ -3.17\!\times\!10^{-15}.
\end{aligned}
\]
while with Algorithm~\ref{alg4} we obtain the eigenvalues
\[
\begin{aligned}
&6.13\!\times\!10^{-15} \pm 7i,\
-0.63,\ -0.43,\ -0.13 \pm 0.99i,\
-0.19 \pm 0.97i,\
-0.49 \pm 0.87i,\\
&-0.44 \pm 0.90i,\
-0.29 \pm 0.95i,\
-0.33 \pm 0.94i,\ 2.32\!\times\!10^{-16},\ -1.84\!\times\!10^{-15}.
\end{aligned}
\]
These results confirm that the purely imaginary pair $ \pm 7i $ is preserved in both reduced systems, while the pair $ \pm 5i $ has been successfully deflated.
}
\end{example}

Our last example demonstrates the positive effect of structure preservation.
\begin{example}\label{ex:rafik}{\rm
Consider randomly 	generated  matrices  $ (M, D, K) $ with $n=10$  $ \rank(D) =5$ and $ \rank(K) =9$. The exact eigenvalues of $P(\lam)$ and the eigenvalues computed by {\tt polyeig(K, D, M)} are given below. \\
	
%
	
{\small	
		\begin{verbatim}
		 exact eigenvalue                                              polyeig(K, D, M)
	 -----------------------                                    ------------------------------
-1.000000000000000e+00 + 0.000000000000000e+00i     -4.713078909233337e-14 + 7.000000000000087e+00i
-1.000000000000000e+00 + 0.000000000000000e+00i     -4.713078909233337e-14 - 7.000000000000087e+00i
-1.500000000000000e+00 - 8.660254037844386e-01i     -1.000000000000015e+00 + 5.916079783099612e+00i
-1.500000000000000e+00 + 8.660254037844386e-01i     -1.000000000000015e+00 - 5.916079783099612e+00i
 0.000000000000000e+00 - 3.000000000000000e+00i     -2.000000000000025e+00 + 4.582575694955889e+00i
 0.000000000000000e+00 + 3.000000000000000e+00i     -2.000000000000025e+00 - 4.582575694955889e+00i
 0.000000000000000e+00 - 4.000000000000000e+00i      3.168255788407744e-14 + 4.000000000000008e+00i
 0.000000000000000e+00 + 4.000000000000000e+00i      3.168255788407744e-14 - 4.000000000000008e+00i
-2.000000000000000e+00 - 4.582575694955840e+00i      9.611363455345153e-15 + 3.464101615137754e+00i
-2.000000000000000e+00 + 4.582575694955840e+00i      9.611363455345153e-15 - 3.464101615137754e+00i
-1.000000000000000e+00 - 5.916079783099616e+00i     -3.893557179285219e-15 + 3.000000000000003e+00i
-1.000000000000000e+00 + 5.916079783099616e+00i     -3.893557179285219e-15 - 3.000000000000003e+00i
 0.000000000000000e+00 - 7.000000000000000e+00i     -1.000000000000003e+00 + 1.732050807568880e+00i
 0.000000000000000e+00 + 7.000000000000000e+00i     -1.000000000000003e+00 - 1.732050807568880e+00i
-1.000000000000000e+00 - 1.732050807568877e+00i     -1.499999999999996e+00 + 8.660254037844435e-01i
-1.000000000000000e+00 + 1.732050807568877e+00i     -1.499999999999996e+00 - 8.660254037844435e-01i
 0.000000000000000e+00 - 3.464101615137754e+00i     -9.999999999999947e-01 + 1.160627480040896e-07i
 0.000000000000000e+00 + 3.464101615137754e+00i     -9.999999999999947e-01 - 1.160627480040896e-07i
 0.000000000000000e+00 + 0.000000000000000e+00i      1.582550664396186e-13 + 0.000000000000000e+00i
 0.000000000000000e+00 + 0.000000000000000e+00i     -1.670817646028280e-13 + 0.000000000000000e+00i
	\end{verbatim}
	}
Note that $\pm 3i$ and $\pm 7i$ are purely imaginary eigenvalues of $P(\lam)$ and $0$ is an eigenvalue of multiplicity $2$.  If we deflate $\pm 3i$, then we have $ \rank((\begin{bmatrix} K-9 M& D\end{bmatrix}) = 9$ which confirms that $ \pm 3 i$ are simple eigenvalues. After deflating $\pm 3 i$, we obtain a reduced order polynomial $\widetilde{P}(\lambda)$ with eigenvalues (computed by {\tt polyeig})
\begin{verbatim}
	-1.125478374718291e-14 + 7.000000000000022e+00i
	-1.125478374718291e-14 - 7.000000000000022e+00i
	-1.000000000000012e+00 + 5.916079783099621e+00i
	-1.000000000000012e+00 - 5.916079783099621e+00i
	-1.999999999999995e+00 + 4.582575694955837e+00i
	-1.999999999999995e+00 - 4.582575694955837e+00i
	 0.000000000000000e+00 + 4.000000000000007e+00i
	 0.000000000000000e+00 - 4.000000000000007e+00i
	-1.585109229061672e-15 + 3.464101615137759e+00i
	-1.585109229061672e-15 - 3.464101615137759e+00i
	-9.999999999999991e-01 + 1.732050807568878e+00i
	-9.999999999999991e-01 - 1.732050807568878e+00i
	-1.670032062499724e-15 + 0.000000000000000e+00i
	-1.499999999999996e+00 + 8.660254037844459e-01i
	-1.499999999999996e+00 - 8.660254037844459e-01i
	-9.999999594690616e-01 + 0.000000000000000e+00i
	-1.000000040530936e+00 + 0.000000000000000e+00i
	 4.289592826761062e-15 + 0.000000000000000e+00i
\end{verbatim}	
Here, the eigenvalues $\pm 3 i$ have been deflated and the other eigenvalues remain almost unaffected. Observe that the double eigenvalue $-1$ of $P(\lam)$  splits into two real eigenvalues very close to $-1$ in the deflated polynomial $\widetilde{P}(\lam)$. In contrast, the double eigenvalue $-1$  of $P(\lam)$  splits into two complex conjugate eigenvalues when eigenvalues of  $P(\lam)$ are computed by {\tt polyeig}. 
}
\end{example}
}
\section{Conclusion}
In this paper, we have proposed a trimmed linearization of a quadratic matrix polynomial
arising from a damped mass-spring system,
in which the eigenvalues of $P$ at infinity, resp., zero and also the other purely imaginary eigenvalues are deflated. Also,
we have presented  results on how one can assess whether the problem is hyperbolic or not, and what is the inertia of the hyperbolic eigenvalue problem.
Finally, we thoroughly describe the movement of purely imaginary eigenvalues under the influence of semidefinite parametric damping matrices.

\subsection*{Acknowledgments:}
The first author thanks his co-authors for inviting him to TU Berlin and University J.J. Strossmayer, Osijek, where parts of this work were prepared. Their generous support and warm hospitality during the visits are gratefully acknowledged.
This paper has been partly supported by the Croatian Science Foundation under the project Conduction, IP-2022-10-5191.

\subsection*{Declaration of competing interest}

The authors declare that they have no known competing financial interests or personal
relationships that could have appeared to influence the work reported in this paper.

\section{Appendix}
In this appendix we present
{\sc Matlab} codes for Algorithms \ref{alg2} and \ref{alg4}.
The corresponding QR based codes for Algorithms \ref{alg1} and \ref{alg3} are similar. We also present two codes for generating specific coefficients $M,D,K$.
	
\begin{verbatim} 





%  R. Alam, V. Mehrmann and N. Truhar     August 10, 2026

 % This program implements SVD based method to construct a trim linearization
 % that deflates 0 eigenvalues of P(s) := s^2 M + s D + K, where K:= GG^*
 % and G is n-by-r matrix with rank(G) =r.


clear

%n= 5;  r= 3; %d = [2, 0, 2, 0, 0]'; w = [ 2, 3, 4 ]'; full deflation

%n= 5;  r= 2; d = [2, 0, 2, 0, 0]'; w = [ 2, 3 ]'; partial deflation

%n= 7; r = 3; d = [2, 0, 2, 0, 0, 0, 0,]'; w = [ 2, 3, 4 ]'; % partial deflation

n= 10; r = 3; d = [2, 0, 2, 0, 0, 0, 0, 0, 0, 0]'; w = [ 2, 3, 4 ]'; % partial deflation


[M, D, G] = genmdk2(n, d, r, w);  % randomly generate M > 0, D >=0 with
                                  % eigenvalues given by d, and G s.t. K= GG^*
                                  % has nonzero eigenvalues given in w.

K= G * G';

%Compute eigenvalues of P(s) = s^2 M + sD + K.

disp('Eigenvalues of P(s) = s^2 M +s D+K')
evpol = polyeig(K, D, M)


% Construct the pencil S(s) = A + s B

 B = [ M, zeros(n, r);
       zeros(r, n), - eye(r)];
 A = [D, G;
      G', zeros(r, r)];

 % compute eigenvalues of S(s)


disp('Eigenvalues of trimmed pencil S(s) = A+sB')

 evpen = eig(A, -B)


 % SVD based method for trim linearization and deflation of 0 eigenvalue

m = rank([D, G]);  % If  m = n then exit as 0 is already deflated.
 if  m ==n;  disp( 'zero eigenvalues already deflated'), return;
 end

[U, S, V] = svd([D, G]);    M  = U'*M*U;

MM = M(1:m, 1:m) - M(1:m, m+1:n) * (M(m+1:n, m+1:n) \ M(1:m, m+1:n)'); % Schur complement of $M_{22}$.
DD  = U(:, 1:m)'* D* U(:, 1:m); % reduced size D.
RR = S(1:m, 1:m) * V(:, 1:m)';  GG = RR(:, n+1:n+r); % last r columns of RR gives reduced size G

% Use MM, DD, GG to construct the timed linearization L(s) = AA + s BB

AA = [ DD, GG;
      (GG)' zeros(r,r)];
BB = [ MM, zeros(m, r);
      zeros(r, m), - eye(r)];


% Compute eigenvalues of the trimmed pencil L(s) = AA + s BB

disp('Eigenvalues of trimmed pencil H(s) = AA+sBB')
evtp = eig(AA, - BB)


% ========== function genmdk2(n, d, r, w) ============
%
% This function generates n-by-n (M , D, K) system with the following
% properties:  (a) M >0, D >= 0  and K = GG^* and G is n-by-r with r=
%  rank(G), eigenvalues of D are given by d and nonzero eigenvalues
% of K= GG^* are given by w.

function [M, D, G] = genmdk2(n, d, r, w)
    % d and w nongenative arrarys of size n and r with r <=n.

     if  n< r;  disp( 'n must be bigger than or equal to r'), return;
     end

    rng(42); A = rand(n); [Q, R] = qr(A);
    %v = randi([5, 20], 1, n)'; %  n positive integer in [5, 20]
    v = randi([1, 10], 1, n)'; %  n positive integer in [.5, 10]
    M = Q * diag(v) * Q'; % positive definite M matrix
    G =  Q(:, 1:r) * diag(sqrt(w));  % GG* positive semidefinite stiffness matrix


    rng(43); A = rand(n); [Q, R] = qr(A);
    D =   Q * diag(d) * Q';  % positive semidefinite damping matrix

end  


% R. Alam, V. Mehrmann and N. Truhar        August 10, 2026

 % This program implements SVD based method to deflate  purely imaginary
 % eigenvalues of P(s) := s^2 M + s D + K, where M>0, D>= 0 and K>=0.


n= 5; d=[ 2 3 0 0 4]'; w = [ 1 3 9 16 0]'; % d and w specifies evalues of D and K

%n= 10; d=[ 2 3 0 0 4, 2,0, 2, 0,0]'; w = [ 1 3 9 16 , 25, 36, 49, 4, 12, 0]';

[ M, D, K, evl] = genmdk(n, d, w); % randomly generates M, D, K
                                   % evl contains exact eigenvalues of
                                   % (M, D, K)

% eigenvalues of (M, D, K) computed by polyeig

evl2 = polyeig(K, D, M);

disp('            Exact eigenvalues of (M, D, K)                 Computed by polyeig(K, D, M)')

disp([evl, evl2])


% check for purely imaginary eigenvalues and then deflate

omega = 3*1i;  %purely imaginary eigenvalues 3i, -3i to be deflated

	m = rank( [K+omega^2 * M, D] );  p = n-m;
     if  p ==0;  disp( 'omega, -omega are not eigenvalues'), return;
     end
         [U, S, V] = svd( [K + omega^2 * M, D] ); % svd
         C = chol(U(:, m+1:n)' *  M * U(:, m+1:n), 'lower'); % Cholesky factorization L*L'
         X2 = U(:, m+1:n) / C' ;   % solves system X2 *  C'  = U(:, m+1:n)  for X2
         [Y, T, W] = svd(M*X2);   L = chol( Y(:, p+1:n)' * M * Y(:, p+1:n), 'lower' );
         X1  = Y(:, p+1:n) / L';  X = [X1, X2];
         MM = X1'*M*X1;  DD = X1'*D*X1;  KK =  X1'*K*X1;

% Reduced system (MM, DD, KK) deflates eigenvalues -\omega and -omega.
% compute eigenvalues of (MM, DD, KK)

disp(' Eigenvalues of (MM, DD, KK) after first round deflation of purely imaginary eigenvalues')
eigenvalue = polyeig(KK, DD, MM)



 % Repeat deflation for the eigenvalues  4i, -4i

 omega = 4*1i; % imaginary eigenvalues 4i and -4i to be deflated
 M= MM; D = DD; K = KK; n=4;

	m = rank( [K+omega^2 * M, D] );  p = n-m;
     if  p ==0;  disp( 'pm i*omega are not eigenvalues'), return;
     end
         [U, S, V] = svd( [K + omega^2 * M, D] ); % svd
         C = chol(U(:, m+1:n)' *  M * U(:, m+1:n), 'lower'); % Cholesky factorization L*L'
         X2 = U(:, m+1:n) / C' ;   % solves system X2 *  C'  = U(:, m+1:n)  for X2
         [Y, T, W] = svd(M*X2);   L = chol( Y(:, p+1:n)' * M * Y(:, p+1:n), 'lower' );
         X1  = Y(:, p+1:n) / L';  X = [X1, X2];
         MM = X1'*M*X1;  DD = X1'*D*X1;  KK =  X1'*K*X1;

 % Reduced system (MM, DD, KK) deflates eigenvalues \omega, - omega.
 % compute eigenvalues of (MM, DD, KK)

      disp(' Eigenvalues of (MM, DD, KK) after second round deflation of purely imaginary eigenvalues')
      eigenvalue2 = polyeig(KK, DD, MM)




      % ===========  function genmdk(n, d, w) =================

% This function generates n-by-n (M , D, K) system which is congruent to
% diag(q_1(s), \ldots, q_n(s)), where q_1(s), ..., q_n(s) are
% given monic quadratic polynomials.

function [M, D, K, evl] = genmdk(n, d, w)
    % d and w nongenative arrarys of size n containing eigenvalues
    % D and K, i.e., D = Udiag(d)U^* and K= Udiag(w) U^*, where U
    %  is a randomly generated nonsingular matrix
    % evl contains exact eigenvalues of (M, D, K) system



    % Generate (M, D, K) system

    rng(42); A = rand(n); [Q, R] = qr(A);
    v = randi([5, 20], 1, n)'; %  n positive integer in [5, 20]
    M = Q*diag(v)*Q'; % positive definite M matrix
    X = Q / diag(sqrt(v)); % then  X' * M *X = I
    D =   X'\ diag(d) / X;  % positive semidefinite damping matrix
    K =  (X'\ diag(w)) / X;  % positive semidefinite stiffness matrix.
                             % with X'*K*X = diag(w)


   % Compute exact eigenvalues of (M,D, K) system as roots of s^2+ sd_i+ w_i

    ev = zeros(2*n, 1);
    for j =1:n
    dc = d(j)^2-4*w(j);  % discriminant
    root1 =  -( d(j) + sqrt(dc))/2;
    root2 =   (-d(j) + sqrt(dc))/2;
    %create odd and even index
     idx1 = 2*j-1;
     idx2 = 2*j;
     ev(idx1) = root1;
     ev(idx2) = root2;
    end

    evl = ev;


end


\end{verbatim}

\end{document}